\documentclass[preprint,3p,times]{elsarticle}

\usepackage{amsmath,amsthm,amsfonts,amssymb,amscd}
\usepackage{graphicx}
\usepackage{caption}
\usepackage{subcaption}
\usepackage{color}

\usepackage{graphicx}
\usepackage{tikz}

\usepackage{algorithm}
\usepackage{algorithmic}
\usepackage{boldfonts}

\usepackage{natbib}
\usepackage{siunitx}

\newtheorem{remark}{Remark}

\newtheorem{form}{Formulation}

\def\ep{\varepsilon}
\def\eps{\varepsilon}
\def\phi{\varphi}

\begin{document}

\begin{frontmatter}

\title{Iterative coupling of flow, geomechanics 
and adaptive  phase-field  fracture
including { level-set crack width approaches}}

\author[ices]{Sanghyun Lee\corref{cor1}}
\ead{shlee@ices.utexas.edu}
\author[ices]{Mary F. Wheeler}
\ead{mfw@ices.utexas.edu}
\author[ricam,ep]{Thomas Wick}
\ead{thomas.wick@polytechnique.edu}

\cortext[cor1]{Corresponding author}

\address[ices]{Center for Subsurface Modeling, The Institute of Computational Engineering and Sciences, The University of Texas at Austin, 201 East 24th Street, Austin TX 78712, USA}
\address[ricam]{Johann Radon Institute for Computational and Applied Mathematics, Austrian Academy of Sciences, 4040 Linz, Austria}
\address[ep]{Centre de Math\'ematiques Appliqu\'ees, \'Ecole Polytechnique,
  91128 Palaiseau, France}

\begin{abstract}
In this work, we present numerical studies of fixed-stress iterative coupling 
for solving flow and geomechanics with propagating fractures in a porous medium. 
Specifically, fracture propagations are described by employing a phase-field approach. 
The extension to fixed-stress splitting to propagating 
phase-field fractures and systematic 
investigation of its properties are important enhancements to existing studies.
Moreover, we provide an accurate computation
of the fracture opening
using level-set approaches and 
a subsequent finite element interpolation of the width. The latter
enters as fracture permeability into 
the pressure diffraction problem which is crucial for fluid filled fractures.
Our developments are substantiated with several 
numerical tests that include comparisons 
of computational cost for iterative coupling and
nonlinear and linear iterations as well as 
convergence studies in space and time.
\end{abstract}


\end{frontmatter}

\section{Introduction}
Iterative coupling has received great importance for coupling 
flow and mechanics in subsurface modeling, environmental and petroleum 
engineering problems
\cite{CaWhiTchel15,KimTchJua11a,KimTchJua11b,MiWaWhe14,MiWhe12,Settari:1998,Settari:2001wf}. 
Recently, the extension of iterative coupling to 
fractured porous media has been of interest 
\cite{Girault2016,MiWheWi14,SiPeKuWiGaWhe14}.
However, reliable and efficient numerical 
methods in coupled poromechanics, including fractures, still pose
computational challenges. The applications include
multiscale and multiphysics 
phenomena such as reservoir deformation,
surface subsidence, well stability, 
sand production, waste deposition, pore collapse, fault
activation, hydraulic fracturing, CO$_2$
sequestration, and hydrocarbon recovery.

On the other hand, 
quasi-static brittle fracture propagation using variational techniques has attracted attention in recent years since the pioneering work in \cite{BourFraMar00,FraMar98}. 
The numerical approach \cite{BourFraMar00} is based 
on Ambrosio-Tortorelli elliptic functionals \cite{AmTo90,AmTo92}. 
Here, discontinuities in the displacement field $\bu$ across the lower-dimensional
crack surface are approximated by an auxiliary function $\varphi$.
This function can be viewed as an indicator function, which
introduces a diffusive transition zone between the broken and the unbroken material. 
This zone has a half bandwidth $\eps$, which is a model
regularization parameter. 
From an application viewpoint, two situations are of interest for given fracture(s):
first, observing the variation of the fracture width (crack opening displacement) and second, change of the fracture length. 
The latter situation is by far more complicated. However, both configurations 
are of importance and variational fracture techniques can be used for both of them.

Fracture evolutions satisfy 
a crack irreversibility constraint such that the resulting 
system can be characterized as a variational inequality. 
Our motivation for employing such a variational approach is that fracture nucleation, propagation, kinking, and
the crack morphology are automatically included in the model.  In addition, explicit remeshing or reconstruction of the crack path is not necessary. 
{ The underlying equations 
are based on principles arising from continuum mechanics  that can be treated with (adaptive) Galerkin finite elements. }
An important  modification of \cite{FraMar98}  towards a thermodynamically-consistent phase-field fracture model has been accomplished in  \cite{MieWelHof10a,MieWelHof10b}. These approaches have been extended  to pressurized fractures in
\cite{MiWheWi15c,MiWheWi15b} that include a decoupled approach 
and a fully-coupled technique accompanied with rigorous analysis. 
Moreover, a free energy functional was established in \cite{MiWheWi15c}.
In these last studies, the 
crack irreversibility constraint has been imposed through penalization.  It it well-known that the energy functional of the basic  displacement/phase-field model is non-convex and constitutes a crucial  aspect in designing efficient and robust methods. 
{  Most approaches for coupling displacement and phase field  are sequential}, e.g., 
\cite{BourFraMar08,Bour07,BuOrSue10,MesBouKhon15}; 
however it is well-known that a monolithic treatment has higher robustness 
(and potentially better efficiency) than sequential coupling. 
Indeed it has been shown in \cite{Vignollet2014,Gerasimov2015}
that for certain phase-field fracture configurations partitioned coupling is more expensive than a monolithic solution. 
Consequently in this paper, we adopt a quasi-monolithic approach 
using an extrapolation in the phase-field variable \cite{HeWheWi15,LeeRebHayWhe_2016}.

{ Recent advances and numerical studies 
for treating multiphysics phase-field fracture
include the following;}
thermal shocks and thermo-elastic-plastic solids \cite{BouMarMauSics14,Miehe2015449,Miehe2015486}, 
elastic gelatin for wing crack formation \cite{LeeRebHayWhe_2016},
pressurized fractures \cite{BourChuYo12,MiWheWi15b,WheWiWo14,WiLeeWhe15,Wi16_dwr_pff}, 
fluid-filled {(i.e., hydraulic)} fractures
\cite{MiWheWi14,LeeWheWi16,Miehe2015186,MieheMauthe2015,Markert2015,Heider2016,Wilson2016264},
proppant-filled fractures \cite{LeeMiWheWi16},
a fractured well-model within a reservoir \cite{WiSiWhe15},
and crack initiations with microseismic probability maps \cite{LeeWheWiSri16}. 
These studies demonstrate that phase-field fracture has great potential to tackle 
practical field problems.

{ 
Addressing multiphysics problems requires careful design of the solution algorithms. 
For the displacement/phase-field subproblem,
we employ a quasi-monolithic approach as previously mentioned, 
but to couple this fractured-mechanics  to flow, we use a splitting approach. 
The latter is more efficient for solvers and for choosing different time scales for both mechanics and flow, respectively. 
In addition, the splitting permits easier extensions to multiphase flow 
including equation of state (EOS) compositional flow.
A successful splitting approach is fixed-stress iterative coupling, which has been applied 
in a series of papers for coupling phase-field fracture/mechanics and flow \cite{MiWheWi14,LeeWheWi16,LeeMiWheWi16,LeeWheWiSri16}.
However a systematic investigation of the performance of this scheme is still missing.
It is the objective of this paper to 
illustrate using benchmarks and mesh refinement studies 
to establish the robustness and efficiency of the fixed-stress algorithm.
These studies are essential for future extensions including efficient three-dimensional practical field problems.
}

In addition, we also focus on a more accurate approximation of the fracture width. 
The authors of \cite{Nguyen2015} recently proposed 
a two-stage level-set approach in which 
first a level-set function is computed with the help of the phase-field  function and in a second step this level-set function is smoothed due to high gradients.
{ 
However, since we only need the level-set function to obtain normal vectors on the fracture boundary, 
we also propose an alternative method which simplifies the  above approach by avoiding the computation of an explicit level-set function but directly using the computed phase-field function.
In addition, we note that these approaches do not derive the width formulation inside the fracture, which is crucial for the  fracture permeability computation for fluid filled fracture propagations.
Thus, here we propose a method to 
compute the crack width values inside the fracture by employing an interpolation based on the values in the diffusive fracture zone. 
}

{ 
To resulting fluid filled fracture propagation framework consists of five equations (four equations
when using phase-field {directly} as a level-set value}) for five {(i.e., four)} unknowns:
vector-valued displacements $\bu$, scalar-valued phase-field $\varphi$, 
pressure $p$, level-set $\varphi_{\textsf{LS}}$, and a finite 
element representation of the fracture width $w$. The first 
problem (namely the displacement/phase-field) is nonlinear
and subject to an inequality constraint-in-time (the crack irreversibility) 
whereas the other three problems are linear.

The outline of this paper is as follows: In Section \ref{sec_math_model}
we recapitulate 
the flow equations in terms of a pressure diffraction problem, 
and the displacement-phase-field system for the mechanics part.
In Section \ref{sec_width}, two equations for computing a level-set function and 
the width are formulated.
In Section \ref{sec_dis} we discuss the discretization of all problems.
In the next Section \ref{sec_fs_splitting} we address 
discretization and the fixed-stress coupling algorithm.
Several numerical examples are presented in Section \ref{sec_tests},
which demonstrate the performance of our algorithmic developments.

\section{Mathematical models for flow and mechanics of porous media and fractures}
\label{sec_math_model}

\subsection{Preliminaries}
\label{sec:prelim}
\begin{figure}[H]
\centering
\includegraphics[scale=0.6]{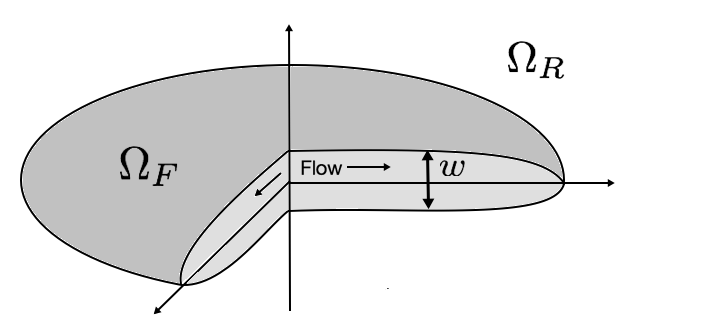}
\caption{Sketch of a penny shaped fracture region ($\mathcal{C}$)  for a three-dimensional setting. The fracture boundary moves in time. }
\label{sketch_fracture_tip}
\end{figure}

\begin{figure}[!h]
\centering
\begin{subfigure}[b]{0.35\textwidth}
\includegraphics[width=\textwidth]{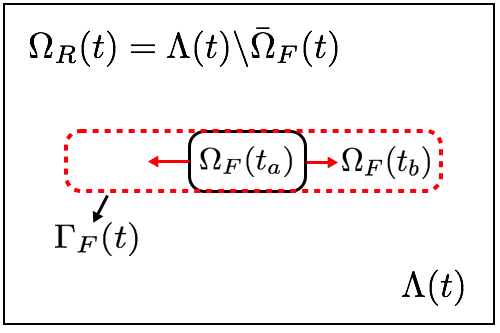}
\caption{}
\end{subfigure}
\hspace*{0.2in}
\begin{subfigure}[b]{0.3\textwidth}
\includegraphics[scale=1.35]{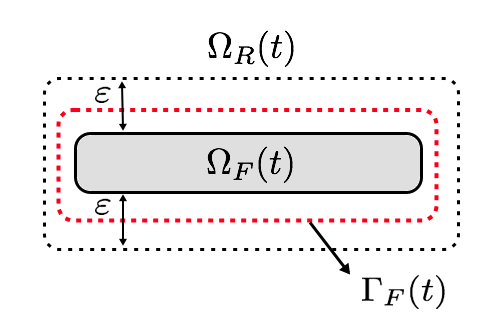}
\caption{}
\end{subfigure}
\caption{Sketch of the (a) time dependent domain for a propagating fracture for two different times $t_a < t_b$ and (b) a fixed fracture with the transition zone where the width is $\varepsilon$. } 
\label{sketch_fracture_diff}
\end{figure}

Let $\Lambda \in \mathbb{R}^d$, $d=2,3$ be a smooth open and 
bounded computational domain with Lipschitz boundary $\partial \Lambda$ and
let $[0,T]$ be the computational time interval {with} $T>0$.
We assume that the crack $\mathcal{C}$ is contained in $\Lambda$. 
The prototype configuration for a horizontal penny shape fracture is given in Figure \ref{sketch_fracture_tip}.
Here, we emphasize that the crack is seen as a thin three-dimensional volume
{$\Omega_F(t)$} using $\varepsilon$ at time $t \in [0,T]$ (see Figure \ref{sketch_fracture_diff}){,}
where the thickness is much larger than the pore size of the porous medium.
The boundary of the fracture 
is denoted by $\Gamma_F(t) := \bar{\Omega}_F(t) \cap \bar{\Omega}_R(t)$, where 
$\Omega_R := \Lambda \setminus \Omega_F$ represents the porous media.

Throughout the paper, we will use the standard notation for Sobolev spaces and their norms. For example, let $E \subseteq \Lambda$, then $\|\cdot\|_{1,E}$ and $|\cdot|_{1,E}$ denote the $H^1(E)$ norm and semi-norm, respectively. 
The $L^2(E)$ inner product is defined as 
$(f, g) := \int_E f g \ d\bx$ for all $v, w \in L^{2}(E)$
with the norm $\| \cdot \|_E$.  
For simplicity, we eliminate the subscripts on the norms if $E = \Lambda$. 
For any vector space $\bX$, $\bX^d$ will denote the vector space of size {$d$}, whose components belong to $\bX$ and $\bX^{d\times d}$ will denote the $d \times d$ matrix whose components belong to $\bX$.

\subsection{A pressure diffraction flow system}
\label{sec_diffraction_pressure}
We now formulate the flow problem in terms 
of a diffraction system \cite{LaSoUr1968}. Specifically, the 
underlying Darcy flow equations have the same structure 
in both the porous medium and the fracture. Using 
varying coefficients and an indicator variable allows 
to distinguish between reservoir flow and fracture flow.

To derive the flow pressure equations for each sub-domain, 
first we consider the two separate mass continuity equations 
for the fluid in the reservoir and the fracture, which we can rewrite as
\begin{align}
\partial_t (\rho_F \phi^\star_F) + \nabla \cdot (\rho_F \bv_F) &= q_F - q_L  
&\text{in } \quad \Omega_F \times (0,T], \label{eqn:density_equation_a} \\
\partial_t (\rho_R \phi^\star_R) + \nabla \cdot (\rho_R \bv_R) &= q_R 
&\text{in } \quad \Omega_R \times (0,T].
\label{eqn:density_equation_b}
\end{align}
Here 
$\rho_{F}$, $\rho_{R}$ are fluid densities,
$q_L$ is a leak-off term (which is assumed to be zero in the following),
$q_F$ and  $q_R$ are source/sink terms for fracture and reservoir, respectively. 

We assume the fluid in the reservoir ($j=R$) and the fracture ($j=F$)  is slightly compressible, 
thus we define the fluid density as
\begin{equation}
\rho_j := \rho_j^0 \exp(c_j (p_j - p_j^0)) 
\approx \rho_j^0 [ 1 + c_j (p_j - p_j^0) ],  \ j \in \{F,R\}
\end{equation}
where 
$p_j: \Omega_j \times [0,T] \rightarrow \mathbb{R}$ is the pressure, 
$p_0$ is the initial pressure at $t=0$,
$\rho_j^0$ is the reference density and $c_j$ is the fluid compressibility. 
In addition, 
$\phi^\star_R$ and $\phi^\star_F$ are the reservoir and fracture fluid
fraction respectively  and we set $\phi^\star_F = 1$ (since the porosity of
the fracture is assumed to be one) and 
\begin{equation}\label{eqn:pressure_main_b}
\phi_R^\star := \phi^\star_0 + \alpha \nabla \cdot \bu + \dfrac{1}{M}(p_R - p_0).
\end{equation}
Here 
$\bu: \Omega \times [0,T] \rightarrow \mathbb{R}^d$ is the solid displacement,  
$\alpha\in [0,1]$ is the Biot coefficient, 
$M>0$ is a given Biot modulus, 
and $\phi^\star_0$ is the initial value.

Next, we describe the flow given by Darcy's law  for the
fracture and for the reservoir , respectively by
\begin{equation}
\bv_j = -\dfrac{K_j}{\eta_j} ( \nabla p_j - \rho_j \bg), \ j \in \{F,R\}
\end{equation}
where $\eta_j$ is the fluid viscosity, $K_j$ is the permeability, 
and $\bg$ is the gravity.

Following the general reservoir approximation with the assumption that 
$c_R$ and $c_F$ are small enough, we use $\rho_R = \rho_R^0$ and 
$\rho_F = \rho_F^0$, and assume $p_0 = 0$,
to rewrite the {equations} \eqref{eqn:density_equation_a}-\eqref{eqn:density_equation_b} by
\begin{align}
&\rho_R^0 \partial_t (\dfrac{1}{M} p_R + \alpha \nabla \cdot \bu ) 
- \nabla \cdot \dfrac{K_R \rho_R^0}{\eta_R} (\nabla p_R - \rho_R^0 \bg) = q_R  
&\text{in } \quad \Omega_R \times (0,T], \label{eqn:pressure_diffraction_a} \\
& \rho_F^0 c_F \partial_t p_F 
- \nabla \cdot \dfrac{K_F \rho_F^0}{\eta_F} (\nabla p_F - \rho_F^0 \bg) = q_F 
&\text{in } \quad \Omega_F \times (0,T] .
\label{eqn:pressure_diffraction_b}
\end{align}
Inside the fracture flow equation, 
the fracture permeability is assumed to be isotropic such that
\begin{equation}
K_F = \frac{1}{12} w(\bu)^2,
\label{eqn:k_f}
\end{equation}
where $w(\bu) = [\bu\cdot \bn]$ denotes the aperture (width) of the fracture,
which means that the jump $[\cdot]$ of normal displacements has to be computed;
corresponding details are provided in Section \ref{sec_width}.
{ For further non-isotropic lubrication laws that have been 
specifically derived for fluid-filled phase-field fractures, we refer to 
\cite{MiWheWi14,LeeMiWheWi16}.}

The system is supplemented with initial and boundary conditions. 
The initial conditions for the pressure diffraction equations
\eqref{eqn:pressure_diffraction_a}-\eqref{eqn:pressure_diffraction_b} are
given by:
{
\begin{align*}
p_F(\bx,0) &=  p_F^0 \quad\text{for all } \bx \in \Omega_F(t=0),\\ 
p_R(\bx,0) &=  p_R^0 \quad\text{for all } \bx \in \Omega_R(t=0), 
\end{align*}
where $p_F^0$ and $p_R^0$ are smooth given pressures.
Also we have 
\[
\phi(\bx,0) = \phi^0 \quad\text{for all } \bx \in \Lambda (t=0), 
\]
where $\phi^0$ is a given smooth initial fracture.
}

We prescribe the boundary and interface conditions for pressure as
\begin{align}
K_{R}(\nabla p_R - \rho_R^0 \bg) \cdot \bn &= 0 
&\text{on} \quad\partial \Lambda \times (0,T], \\
[p] &= 0  &\text{on }  \hspace{0.23in} \Gamma_F \times (0,T],  \\
\dfrac{K_{R}\rho_R^0}{\eta_R}(\nabla p_R - \rho_R^0 \bg)  \cdot \bn &= \dfrac{K_{F}\rho_F^0}{\eta_F}(\nabla p_F - \rho_F^0 \bg)  \cdot \bn &\text{on }  \hspace{0.23in} \Gamma_F \times (0,T],  
\label{eqn:boundary_cond}
\end{align}
where $\bn$ is the outward pointing unit normal on $\Gamma_F$ or $\partial \Lambda$.

In order to finalize our derivation we perform two steps; 
first, we introduce the coefficients with indicator functions to
combine equations \eqref{eqn:pressure_diffraction_a} and \eqref{eqn:pressure_diffraction_b}, 
secondly, we formulate the weak form in terms of a pressure diffraction system. 
The weak formulation reads: 
\begin{form}
\label{form_pressure_weak}
Find $p(\cdot, t)  \in V_p = H^1(\Lambda)$ for almost all times $t\in (0,T]$ such that,
 \begin{equation}
      \rho_0(\partial_t {\phi^\star},v)  +
      ({\rho_0 K_{eff}}(\nabla p - \rho_0 \bg ), \nabla v) - (q, v)
       = 0, \quad
      \forall v \in V_p
      \label{eqn:pressure_diffrac}
  \end{equation}
where the coefficient functions are defined as
\begin{align}
\rho_0 &=\chi_{\Omega_R} \rho_R^0 + \chi_{\Omega_F} \rho_F^0,\\
{\phi^\star} &= {\phi^\star}(\cdot,t) := 
\chi_{\Omega_R}  \Bigl( \dfrac{1}{M}p_R + \alpha \nabla \cdot \bu  \Bigr)
+ \chi_{\Omega_F} ( c_F p_F),\\
{q} &= q(\cdot,t) := \chi_{\Omega_R}{q_R} + \chi_{\Omega_F}  q_F,\\   
K_{eff} &:=\chi_{\Omega_R} \dfrac{K_R}{\eta_R} + \chi_{\Omega_F} \dfrac{K_F}{\eta_F}, 
\end{align}
where $\chi_{\Omega_R}=1$  and $\chi_{\Omega_F}=0$ in $\Omega_R$, and 
$\chi_{\Omega_R}=0$ and 
$\chi_{\Omega_F}=1$ in $\Omega_F$ .
\end{form}

\subsection{Geomechanics and phase-field fracture equations}
\label{sec_pff}

The displacement of the solid and diffusive flow in a non-fractured porous medium are modeled in 
$\Omega_R$ by the classical quasi-static elliptic-parabolic 
Biot system for a 
porous solid saturated with a slightly compressible viscous fluid.
The constitutive equation for the Cauchy stress tensor is given as 
\begin{equation}
\sigma^{por}(\bu, p)  -\sigma_0 = \sigma(\bu) - \alpha (p - p_0) I,
\label{eqn:model:cauchystress}
\end{equation}
where $I$ is the identity tensor and 
$\sigma_0$ is the initial stress value.
The effective linear elastic stress tensor is 
\begin{equation}
\sigma := \sigma(\bu) = \lambda(\nabla \cdot \bu) I + 2 G e(\bu),
\end{equation}
where $\lambda, G >0$ are the Lam\'e coefficients.
The linear elastic strain tensor is given as
$e(\bu) := \frac{1}{2}(\nabla \bu + \nabla \bu^T)$.
Then the balance of linear momentum in the solid reads 
\begin{equation}
-\nabla \cdot \sigma^{por}(\bu, p) = \rho_s \bg  
\quad \text{in }  \Omega_R \times (0,T] ,
\label{eqn:model:linearmomen}
\end{equation}
where $\rho_s$ is the density of the solid.
We prescribe homogeneous Dirichlet boundary conditions on $\partial\Lambda$ for the displacement $\bu$.

In the following, we describe our fracture approach in a porous medium using
the previous setup.
Modeling fractures with a phase-field approach in $\Lambda$ is formulated 
with the help of an elliptic (Ambrosio-Tortorelli) functional \cite{AmTo90,AmTo92}
and a variational setting,
which has been first proposed for linear elasticity in \cite{FraMar98,BourFraMar00}. 

We now recapitulate the essential elements for 
a phase-field model for pressurized and fluid filled fractures in porous media, 
which has been modeled in \cite{MiWheWi15c,MiWheWi15b} including rigorous analysis.
Two unknown solution variables are sought, namely
vector-valued displacements $\bu(\cdot, t)$
and a smoothed scalar-valued indicator phase-field function $\varphi(\cdot, t)$.
Here $\varphi = 0$ denotes the crack region ($\Omega_F$)
and $\varphi = 1$ characterizes the unbroken 
material ($\Omega_R$). The intermediate values constitute a smooth transition zone 
dependent on a regularization parameter $\ep > 0$.

The physics of the underlying problem requires 
a crack irreversibility condition 
that is an inequality condition in time:
\begin{align}
\label{eq_crack_irre}
\partial_t \varphi \leq 0.
\end{align}
Consequently, modeling of fracture evolution problems leads to a variational 
inequality system, that is always, due to this constraint,
quasi-stationary or time-dependent.

The resulting variational formulation 
is stated in an incremental (i.e., time-discretized) formulation
in which the continuous irreversibility constraint is approximated by
\[
\varphi \leq \varphi^{old}.
\]
Here, $\varphi^{old}$ will later denote the previous time step solution and
$\varphi$ the current solution. Let the function spaces be given by 
$V:=H^1_0(\Lambda), W:=H^1(\Lambda)$ and 
\[
W_{in}:=\{w\in H^1(\Lambda) |\, w\leq \varphi^{old} \leq 1 \text{ a.e. on } \Lambda\}.
\] 
We note that the phase field function is subject to homogeneous Neumann conditions on $\partial\Lambda$.
The Euler-Lagrange system for pressurized phase-field fracture reads \cite{MiWheWi15b}:
\begin{form}
\label{form_1a}
Let $p\in H^1(\Lambda)$ be given. Find $\{\bu,\varphi\} \in V \times W$ such that 
\begin{equation}\label{E22ATW}
\begin{aligned}
  &\Bigl(\big( (1-\kappa) {\varphi}^2  +\kappa \big)\;\sigma^+(\bu), e( {\bw}
  )\Bigr)  
+ (\sigma^-(\bu), e(\bw)) \\
&\quad -(\alpha - 1)({\varphi}^{2} p, \mbox{div }  {\bw}) 
+ ({\varphi}^{2} \nabla p,  {\bw}) 
=0 \quad \forall \bw\in V ,
\end{aligned}
\end{equation}
and
\begin{equation} \label{E33ATW}
\begin{aligned}
 &(1-\kappa) ({\varphi} \;\sigma^+(\bu):e( \bu), \psi {-\varphi}) 
-  2(\alpha - 1) ({\varphi}\;  p\; \mbox{div }  \bu,\psi{-\varphi})
+ 2\, ({\varphi} \nabla p\;  {u},\psi) 
\\
&+  G_c  \Bigl( -\frac{1}{\ep} (1-\varphi,\psi{-\varphi}) + \ep (\nabla
\varphi, \nabla (\psi - {\varphi}))   \Bigr)  \geq  0
\quad \forall \psi \in W_{in}\cap L^{\infty}(\Lambda).
\end{aligned}
\end{equation}
\end{form}
Here, $G_c$ is the critical energy release rate
and $\kappa$ is a very small positive regularization parameter ($\kappa \approx 0$) for the 
elastic energy {(in some cases, see for instance \cite{BoVeScoHuLa12}, $\kappa
= 0$ even works)}.
Physically, $\kappa$ represents the residual stiffness of the
material. Consequently, since  
\[
\big( (1-\kappa) {\varphi}^2  +\kappa \big) \; \rightarrow \kappa\quad\text{for
}
\varphi \to 0,
\]
the material stiffness decreases while approaching the fracture zone.
Regarding the stress tensor split, we follow \cite{Amor_2009} in which
the stress tensor is additively decomposed into a tensile part $\sigma^+(\bu)$ 
and a compressive part $\sigma^-(\bu)$ by: 
\begin{align}
\sigma^+(\bu) &:= (\dfrac{2}{n} G + \lambda) tr^+(e(\bu)) I + 2 G (e(\bu) - \dfrac{1}{n} tr(e(\bu)) I ), \\
\sigma^-(\bu) &:= (\dfrac{2}{n} G + \lambda) tr^-(e(\bu)) I,
\label{eqn:stress_split}
\end{align}
where $n$ is the dimension (2 or 3) and 
\begin{equation}
tr^+(e(\bu)) = \max(tr(e(\bu)),0), \quad tr^-(e(\bu)) = tr(e(\bu)) - tr^+(e(\bu)).
\end{equation}
We emphasize that the energy degradation only acts on the tensile part.

\section{The fracture width computation using a level-set approach}
\label{sec_width}
A crucial issue in fluid-filled fractures is the fracture width $w:=w(\bu)$ computation
since this enters as fracture permeability values \eqref{eqn:k_f} into the pressure diffraction 
problem.
Previously, the fracture width was approximated by given point values of displacements in  \cite{LeeWheWi16,LeeMiWheWi16} 
and a method for a single fracture constructing an additional displacement field was studied in \cite{VeBo13}.
In this section, we introduce a computational method to compute the crack width more robustly, especially inside the fracture.

{ The first challenge is to compute the normal vector of the fracture interface. 
Here the method is
inspired by a recent idea proposed in \cite{Nguyen2015} by introducing
an {(explicit)} level-set function $\varphi_{\textsf{LS}}$ for the fracture. 
However, we also note that employing a level-set function to compute the normal vectors of iso-surfaces has been used for many different applications (for example, see  \cite{hysing2008numerical,BonGueLee2016} and references cited therein). 
}
The initial step follows 
a standard procedure in level-set methods.
Let $\Gamma_F$ be the fracture boundary. We now define 
$\Gamma_F$ as the zero levels-set of a function 
$\varphi_{\textsf{LS}}$ such that
\begin{align*}
\varphi_{\textsf{LS}} &> 0, \quad \bx\in\Omega_R,\\
\varphi_{\textsf{LS}} &< 0, \quad \bx\in\Omega_F,\\
\varphi_{\textsf{LS}} &= 0, \quad \bx\in\Gamma_F,
\end{align*}
where $\Gamma_F := \{ \bx \in \Lambda \ | \ \phi(\bx,t) = { C_{LS}}  \}$, 
$\Omega_R := \{ \bx \in \Lambda \ | \ \phi(\bx,t) > C_{LS}  \}$
and  $\Omega_F := \{ \bx \in \Lambda \ | \ \phi(\bx,t) < C_{LS}  \}$.
{ 
Here $C_{LS} \in (0,1)$ is a constant that we have to choose to define the fracture boundary $\Gamma_F$, since
the phase field approach involves the diffusion zone between $\Omega_R$ and $\Omega_F$ with length $\varepsilon$ (see Figure \ref{sketch_fracture_diff}).
However, since $\varepsilon$ is very small and the dependency of the choice of $C_{LS}$ is minimal (e.g \cite{LeeWheWi16}), 
we set  $C_{LS} = 0.1$ throughout this paper for simplicity.
{In the following,} we propose two different methods to compute
$\varphi_\textsf{LS}$, {where one (Formulation \ref{form:level_set}) is a similar technique} as shown in \cite{Nguyen2015}.
}

\begin{figure}[!h]
\centering
\begin{subfigure}[b]{0.25\textwidth}
{\includegraphics[width=\textwidth]{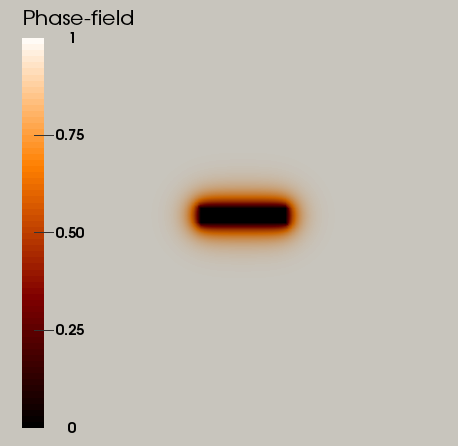}}
\caption{Phase field $\varphi$.}
\label{fig:sec3_fracture_a}
\end{subfigure}
\hspace*{0.1in}
\begin{subfigure}[b]{0.25\textwidth}
{\includegraphics[width=\textwidth]{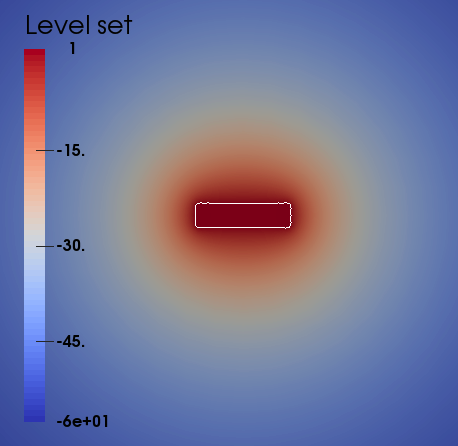}}
\caption{Level-set $\varphi_\textsf{LS}$}
\label{fig:sec3_fracture_c}
\end{subfigure}
\hspace*{0.1in}
\begin{subfigure}[b]{0.25\textwidth}
{\includegraphics[width=\textwidth]{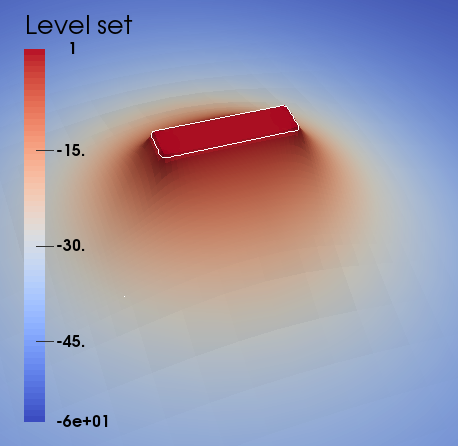}}
\caption{Level-set $\varphi_\textsf{LS}$}
\label{fig:sec3_fracture_b}
\end{subfigure}
\caption{For a given fracture in the middle of domain, 
(a) the phase field values and  in (b)-(c) the level-set values with zero
  level-set ($\varphi_\textsf{LS}=0)$ are illustrated. In (c), a three-dimensional surface
  plot is shown.}
\label{fig:sec3_fracture}
\end{figure}

\begin{form}[Level-set values {obtained by computing an additional problem}]
\label{form:level_set}
Find $\varphi_{\textsf{LS}}$ such that
\begin{align}
-\Delta\varphi_{\textsf{LS}} &= f_{\textsf{LS}}(\cdot,t)\quad\text{in } \Lambda,\\
\varphi_{\textsf{LS}} &= 0 \quad\text{on } \Gamma_F,\\
\partial_{\bn} \varphi_{\textsf{LS}} &= 0 \quad\text{on } \partial\Lambda,
\label{eqn:levelset_system}
\end{align}
where 
\[
f_{\textsf{LS}}(\cdot,t) = \chi(\cdot,t) f_1 + (1-\chi(\cdot,t)) f_2
\]
with $\chi(\cdot,t) = 0$ for $\varphi(\bx,t) < C_{LS}$ and $\chi(\cdot,t) = 1$ 
otherwise. For simplicity, we set $f_1= -10$ and $f_2= 10$.
\end{form}
Figure \ref{fig:sec3_fracture} illustrates the result for the level-set formulation with a simple fracture in the middle of the domain.

{ The other alternative approach differs from Formulation \ref{form:level_set} is to directly using the phase-field values by shifting them but without solving an additional problem:}
{
\begin{form}[Level-set values obtained from phase-field]
\label{form_level_set_by_PFF}
As second alternative approach, the level-set values are 
immediately obtained from the phase-field by:
\[
\varphi_{\textsf{LS}} = \varphi - C_{LS}.
\]
\end{form}
}

{
\begin{remark}[Regularity properties of $\varphi_{\textsf{LS}}$]
Exemplarily, we refer to \cite[section 4.5]{hysing2008numerical}
  for discussions regarding to the regularity properties of
  $\varphi_{\textsf{LS}}$ and to improve the accuracy for computing the gradients across the $\Gamma_F$.
\end{remark}
}

Next, with the computed level-set value $\varphi_{\textsf{LS}}$ {(either
obtained from Formulation \ref{form:level_set} or \ref{form_level_set_by_PFF}),} we obtain the outward
normal vector for a given 
level-set fracture boundary by the following procedure:
\begin{form}[Computing the width with the normal vector on the fracture boundary]
Under the assumption that $\bu^+\cdot \bn = -\bu^-\cdot \bn$ (symmetric displacements
at the fracture boundary), we compute the width locally in each quadrature point
\begin{equation}
w_D := 2\bu\cdot \bn_F = -2 \bu\cdot \frac{\nabla\varphi_{\textsf{LS}}}{\|\nabla\varphi_{\textsf{LS}}\|} \quad\text{on }\Gamma_F.
\label{eqn:new_width}
\end{equation}
Here we assume that {$\bu^+\cdot\bn = -\bu^- \cdot\bn$}, which has been justified for tensile stresses and homogeneous isotropic media. 
These results are compatible for fluid filled fracture 
as we observe in our computational results, e.g., see Figure \ref{ex_1_fig_2}, in comparison to 
just taking two times the $\bu_y$ displacements.
\end{form}

{ The second challenge is to compute the width values inside the crack, which is required for the permeability in fluid filled fracture propagation, see \eqref{eqn:k_f}.
Here we propose a method, Formulation \ref{form:width}, to 
compute the crack width values inside the fracture by employing an interpolation based on the values on $\Gamma_F$. 
}
\begin{form}[Crack width interpolation inside the fracture]
\label{form:width}
We solve the following width-problem: Find $w\in H^1$ such that
\begin{align}
-\Delta w &= g \quad \text{in }\Lambda, \nonumber \\
w &= w_D \quad\text{on } \Gamma_F, \label{eqn:form3}  \\
w &= 0 \quad\text{on } \partial\Lambda. \nonumber 
\end{align}
Here {{$g(\bx) = \beta \|w\|_{L^{\infty}(\Lambda)}$, where $\beta \approx 100$,}} in order to 
obtain a smooth parabola-type width-profile in the fracture. 
{ We note that $\beta$ is problem-dependent and heuristically chosen.}
In the case of multiple fractures, say $m$ fractures, we determine 
a locally highest width, where
$\|w\|_{L^{\infty}(\Lambda_{l})}, \ l \in \{1,\cdots, m\}$ is defined on the local region near the fracture such that $\Lambda := \Lambda_{1} \cup \Lambda_{2} \cup \cdots \cup \Lambda_{m}$. 
\end{form}

\begin{figure}[H]
\centering
\begin{subfigure}[b]{0.225\textwidth}
{\includegraphics[width=\textwidth]{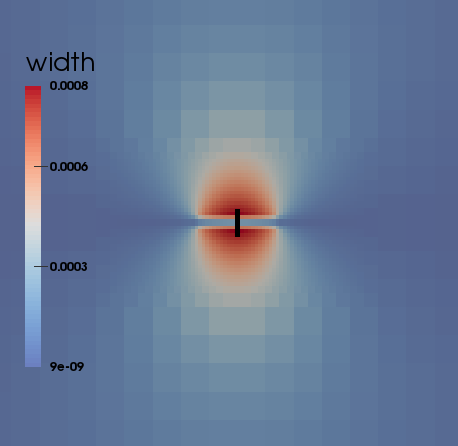}}
\caption{$w_D$ by \eqref{eqn:new_width}.}
\label{fig:sec3_fracture_a}
\end{subfigure}
\hspace*{0.1in}
\begin{subfigure}[b]{0.225\textwidth}
{\includegraphics[width=\textwidth]{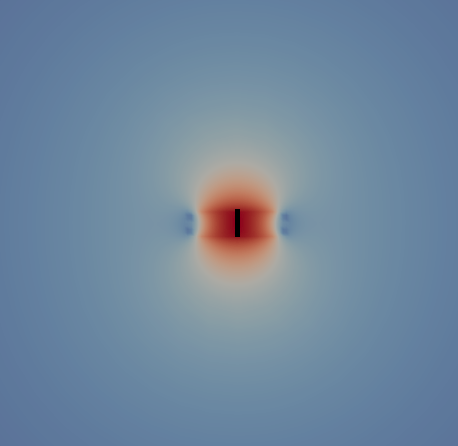}}
\caption{$w$ by \eqref{eqn:form3}. }
\label{fig:sec3_fracture_b}
\end{subfigure}
\hspace*{0.1in}
\begin{subfigure}[b]{0.32\textwidth}
{\includegraphics[width=\textwidth]{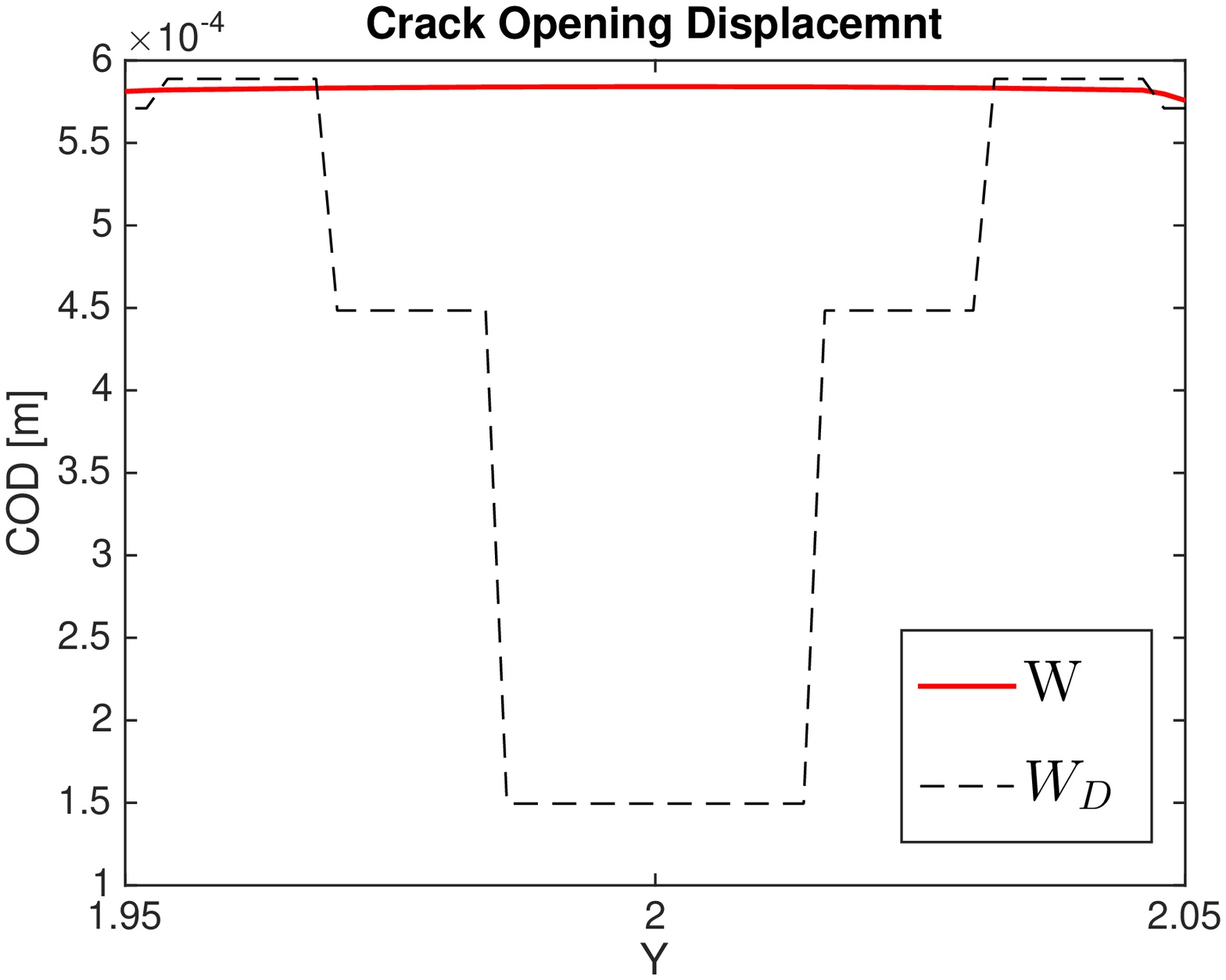}}
\caption{Comparison.}
\label{fig:sec3_fracture_c}
\end{subfigure}
\caption{In (a), we visualize a fracture width evaluation by employing \eqref{eqn:new_width} 
for the phase field fracture and the level-set given in Figure \ref{fig:sec3_fracture}.
The fracture width is computed correctly up to the fracture boundary. 
However, in the very inner of the fracture $\varphi \simeq 0$, all values are nearly zero. 
Then, (b) we compute the width in each
quadrature point and finally solve a width-problem \eqref{eqn:form3} in order to interpolate 
the fracture-boundary width values inside the fracture. (c) shows the comparison between (a) and (b) over the middle line $y-$direction illustrated in the fracture. }
\label{fig:sec3_fracture_width}
\end{figure}

Figure \ref{fig:sec3_fracture_width} illustrates the result for {Formulations \ref{form:level_set}-\ref{form:width}} with crack width in the fracture. Finally we can approximate accurate crack width values in the fracture.

\begin{remark}
{ Using Formulation \ref{form:level_set}} for computing the width is accompanied by the cost that we need 
to solve two additional problems.  
However, we emphasize that in our algorithm these two subproblems are scalar-valued, 
linear, and elliptic 
and therefore much cheaper to compute in comparison 
to the other subproblems.
The advantage of this procedure being that we 
have an accurate width computation as well as a representation  of a global finite element function that can be easily accessed in the program.
In addition, the above computation can be employed 
for multiple non-planar fractures. 
{
To further reduce the computational cost Formulation
\ref{form:level_set} {and can be replaced by Formulation 
\ref{form_level_set_by_PFF}}.
}
\end{remark}

\section{Discretization of all sub-systems}
\label{sec_dis}
In this section, we first address discretization of the 
pressure flow system and then we consider 
the displacement-phase-field system.
Finally, we provide the variational formulations
of the discrete level-set and width problems.
We consider a mesh family $\{\mathcal{T}_h \}_{h>0}$, which is assumed to be shape regular in the sense of Ciarlet, and we 
assume that each mesh $\mathcal{T}_h$ is a subdivision of $\bar{\Lambda}$ made of disjoint elements $\mathcal{K}$, i.e., squares when  $d=2$ or cubes when $d=3$. 
Each subdivision is assumed to exactly approximate the computational domain, thus $\bar{\Lambda} = \cup_{K\in\mathcal{T}_h} \mathcal{K}$. 
The diameter of an element $\mathcal{K}\in \mathcal{T}_h$ is denoted by $h$
and we denote $h_{\min}$ for the minimum. For any integer $k \geq 1$ and any
$\mathcal{K} \in \mathcal{T}_h$, we denote by $\mathbb{Q}^k(\mathcal{K})$ the
space of scalar-valued multivariate polynomials over $\mathcal{K}$ of partial
degree of at most $k$. The vector-valued counterpart of $\mathbb{Q}^k(\mathcal{K})$ is denoted 
$\pmb{\mathbb{Q}}^k(\mathcal{K})$.
We define a partition of the time interval 
$0 =: t^0 < t^1 < \cdots < t^N :=T$ and denote the time step size by $\delta t := t^n - t^{n-1}$.

\subsection{Decomposing the domain $\Lambda$ into $\Omega_R$ and $\Omega_F$}

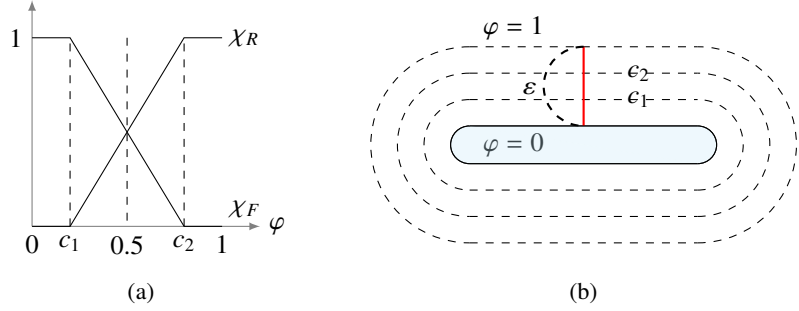
\begin{figure}[h]
  \centering
  \begin{subfigure}[b]{0.225\textwidth}
   \begin{tikzpicture}[scale = 2.5]
    \coordinate (End)   at (1,1);
    \coordinate (XAxisMin) at (0,0);
    \coordinate (XAxisMax) at (1,0);
    \coordinate (YAxisMin) at (0,0);
    \coordinate (YAxisMax) at (0,1);
    \coordinate (x1) at   (0.2,0);
    \coordinate (x1-1) at (0.2,1);
    \coordinate (x2) at   (0.8,0.);
    \coordinate (x2-1) at (0.8,1.);
    
    \coordinate (middle1) at (0.5,1.);
    \coordinate (middle2) at  (0.5,0.);
    
    \draw [thin, gray,-latex] (XAxisMin) -- (1.2,0);
    \draw [thin, gray,-latex] (YAxisMin) -- (0,1.2);
    \draw [thin, dashed] (x1) -- (x1-1);        
    \draw [thin, dashed] (x2) -- (x2-1);
    
    \draw [thin, dashed] (middle1) -- (middle2);        

    \draw [] (middle2) node [below] {$0.5$};
        
    \draw [] (1.2,0) node [right] {$\phi$};
    
    \draw [] (x1) node [below] {$c_1$};
    \draw [] (x2) node [below] {$c_2$};   
    
    \draw [] (YAxisMin) node [below] {$0$};
    \draw [] (YAxisMax) node [left] {$1$};    
    
    \draw [] (XAxisMin) -- (x1) -- (x2-1) -- (End);
    \draw [] (End) node [right] {$\chi_{R}$};
    \draw [] (YAxisMax) -- (x1-1) -- (x2) --  (XAxisMax);
    \draw [] (XAxisMax) node [above right] {$\chi_{F}$};
    
    \draw [] (XAxisMax) node [below] {$1$};

  \end{tikzpicture}
   \caption{}
   \end{subfigure}
   \hspace*{0.2in}
\begin{subfigure}[b]{0.4\textwidth}
\centering
\begin{tikzpicture}

\draw[] (-0.5,0) -- (2.5,0);
\draw[] (-0.5,0.5) -- (2.5,0.5);
\draw (2.5,0) arc (-90:90:0.25) ;
\draw (-0.5,0.5) arc (90:270:0.25) ;

\draw[dashed]  (-0.5,-0.35) -- (2.5,-0.35);
\draw[dashed]  (-0.5,0.85) -- (2.5,0.85);
\draw[dashed]  (2.5,-0.35) arc (-90:90:0.6) ;
\draw[dashed]  (-0.5,0.85) arc (90:270:0.6) ;

\draw[dashed] (-0.5,-0.7) -- (2.5,-0.7);
\draw[dashed] (-0.5,1.2) -- (2.5,1.2);
\draw[dashed]  (2.5,-0.7) arc (-90:90:0.95) ;
\draw[dashed]  (-0.5,1.2) arc (90:270:0.95) ;

\draw[dashed] (-0.5,-1.05) -- (2.5,-1.05);
\draw[dashed] (-0.5,1.55) -- (2.5,1.55);
\draw[dashed]  (2.5,-1.05) arc (-90:90:1.3) ;
\draw[dashed]  (-0.5,1.55) arc (90:270:1.3) ;

\draw[color=red,thick] (1,0.5) -- (1.,1.55);
\draw[dashed,thick] (1,1.55) arc (90:270:0.525) ;

\draw[] (0.5,1.) node[left]  {$\varepsilon$};

\draw[] (0.6,0.25) node[left]  {$\phi=0$};
\draw[] (0.6,1.8) node[left]  {$\phi=1$};
\draw[] (0.6,-1.2) node[left]  {};

\draw[] (2.,0.85) node[left] {$c_1$};
\draw[] (2.,1.2) node[left]  {$c_2$};

\filldraw[fill=cyan!20, fill opacity=.3] (-0.5,0.5) -- (2.5,0.5) arc (90:-90:0.25) -- (2.5,0);
\filldraw[fill=cyan!20, fill opacity=.3] (2.5,0.) -- (-0.5, 0.) arc (270:90:0.25) -- (-0.5,0.5);
\end{tikzpicture}
\caption{}
\end{subfigure}   
\caption{(a) The linear indicator functions $\chi_F$ and $\chi_R$ illustrated with adjustable constants $c_1$ and $c_2$. (b) We consider as the fracture zone if $\phi \leq c_1$ and as the reservoir zone if $\phi \geq c_2$.}
\label{fig:chi_R_and_chi_F}
\end{figure}

We define the fracture domain $\Omega_F$ and the reservoir domain $\Omega_R$ by introducing two linear indicator functions $\chi_F$ and $\chi_R$ for 
the two different sub-domains; they satisfy
\begin{align}
\chi_R (\cdot, \phi) &:= \chi_R(\bx, t, \phi) = 1 
\quad \text{in }  \quad \Omega_R(t),
\quad
\text{and } \quad \chi_R(\cdot, \phi) =  0 
\quad \text{in } \quad  \Omega_F(t), \label{eq_sec_4_chi_a} \\
\chi_F  (\cdot, \phi)&:= \chi_F(\bx, t, \phi) = 1 
\quad \text{in } \quad  \Omega_F(t),
\quad 
\text{and } \quad \chi_F(\cdot, \phi) =  0 
\quad \text{in } \quad  \Omega_R(t). 
\label{eq_sec_4_chi_b}
\end{align}
Thus $\chi_F(\cdot, \phi)$ is zero {in} the reservoir domain and $\chi_R(\cdot, \phi)$ is zero {in} the fracture domain. In the diffusive zone, the linear functions are {defined} as 
\begin{equation}
\label{eq_sec_4_chi_c}
\chi_F(\cdot, \phi) =  -\dfrac{( \phi - c_2 )}{(c_2 - c_1)}
\quad
\text{and} 
\quad
\chi_R(\cdot, \phi) =   \dfrac{( \phi  -  c_1  )}{( c_2 - c_1  )}.
\end{equation}
Thus $\chi_R(\cdot, \phi) =0$ and $\chi_F(\cdot, \phi) =1 $ if $\phi(\bx,t) \leq c_1$, and 
$\chi_R(\cdot, \phi) =1$ and $\chi_F(\cdot, \phi) =0 $ if $\phi(\bx,t) \geq c_2$,
where
$c_1 := 0.5 - c_x$ and  $c_2 = 0.5 + c_x$. 
For simplicity we set $c_x= 0.1$.

\subsection{Temporal and spatial discretization of the pressure diffraction
equation}
The space approximation $P$ of the pressure function $p(\bx,t)$ is approximated by using continuous piecewise polynomials given in the finite element space,
\begin{equation}
\mathbb{W}(\mathcal{T}) := \{ W \in C^0(\bar{\Lambda};\mathbb{R}) \ | \ W|_{\mathcal{K}} \in \mathbb{Q}^1(\mathcal{K}), \forall \mathcal{K} \in \mathcal{T} \}.
\end{equation}
Note that we employ enriched Galerkin approximation \cite{LeeEG16} when we couple with transport problem which requires local and global conservation for flux \cite{LeeMiWheWi16}.

Assuming that the displacement field $\bu$
and the phase field $\varphi$ 
are known, the Galerkin approximation of  \eqref{eqn:pressure_diffraction_a}-\eqref{eqn:pressure_diffraction_b} is formulated as follows. 
Given $P(\bx,0) = P^0$ where $P^0$ is an approximation of the initial condition $p^0$, find $P \in C^1([0,T];\mathbb{W}(\mathcal{T}))$ such that 
\begin{equation}
\chi_R(\cdot, \phi) \Biggl(  \int_{\Lambda}\rho_R^0 \partial_t (\dfrac{1}{M} P + \alpha \nabla \cdot \bu )\omega \ d\bx
+ \int_{\Lambda}\dfrac{K_R \rho_R^0}{\eta_R} (\nabla P - \rho_R^0 \bg) \nabla \omega \ d\bx
=\int_{\Lambda} q_R  \omega  \ d\bx \Biggr),
\quad \forall \omega \in \mathbb{W}(\mathcal{T}),
\end{equation}

\begin{equation}
\chi_F(\cdot, \phi) \Biggl( \int_{\Lambda} \rho_F^0 c_F \partial_t P \omega \ d\bx
+ \int_{\Lambda} \dfrac{K_F \rho_F^0}{\eta_F} (\nabla P - \rho_F^0 \bg)  \nabla \omega \ d\bx= 
\int_{\Lambda} (q_F  - q_L)\omega \ d\bx \Biggr), 
\quad \forall \omega \in \mathbb{W}(\mathcal{T}).
\end{equation}

We denote the approximation of $P(\bx,t^n)$, $0 \leq n \leq N$ by $P^n$, and 
assume $\bu(t^{n+1})$ and $\phi(t^{n+1})$ are given values at time $t^{n+1}$.
Then, the time stepping proceeds as follows: 
Given $P^n$, compute $P^{n+1} \in \mathbb{W}(\mathcal{T})$ so that

\begin{multline}
A_{PR}(P^{n+1})(\omega) := 
\chi_R(\cdot, \phi(t^{n+1})) \Biggl( \int_{\Lambda}\rho_R^0 \Big(\dfrac{1}{M} \Big( \dfrac{P^{n+1} - P^n}{\delta t}  \Big)  
+ \alpha { \Big(\frac{\nabla \cdot \bu^{n+1} - \nabla\cdot \bu^n}{\delta t} \Big)} \Big)  \omega \ d\bx  \\
+ \int_{\Lambda}\dfrac{K_R \rho_R^0}{\eta_R} (\nabla P^{n+1} - \rho_R^0 \bg) \nabla \omega \ d\bx 
- \int_{\Lambda} q_R  \omega  \ d\bx  \Biggr)
\quad \forall \omega \in \mathbb{W}(\mathcal{T})
\end{multline}

\begin{multline}
A_{PF}(P^{n+1})(\omega) :=
\chi_F(\cdot, \phi(t^{n+1})) \Biggl( \int_{\Lambda} \rho_F^0 c_F \Big( \dfrac{P^{n+1} - P^n}{\delta t}  \Big)  \omega \ d\bx
+ \int_{\Lambda} \dfrac{K_F \rho_F^0}{\eta_F} (\nabla P^{n+1} - \rho_F^0 \bg)  \nabla \omega \ d\bx  \\
- \int_{\Lambda} (q_F  - q_L)\omega \ d\bx \Biggr), 
\quad \forall \omega \in \mathbb{W}(\mathcal{T}).
\end{multline}

\begin{form}
\label{form_1}
Find $P^{n+1} \in \mathbb{W}(\mathcal{T})$ for  $t^{n+1}, n=0,1,2,\ldots$ such that
\begin{equation}
A_P(P^{n+1})(\omega) = A_{PR}(P^{n+1})(\omega) + A_{PF}(P^{n+1})(\omega) = 0 \quad\forall\omega\in\mathbb{W}(\mathcal{T}).
\end{equation}
\end{form}

\begin{remark}
\label{remark_permea_at_tip}
To avoid
a singular behavior at the fracture tip in modeling, computations and loss of regularity,
we determine a permeability-viscosity ratio $K_{eff}$ by
interpolation; a so-called cake region
\cite{Kov10} that
is determined by the phase-field variable. 
Our definition of the
cake region is defined in \eqref{eq_sec_4_chi_a} - \eqref{eq_sec_4_chi_c}.
Specifically, outside the cake region, we use
in the reservoir $K_{eff} = K_R/\eta_R$ and in the fracture
$K_{eff} = K_F/\eta_F$.
The resulting interpolated
permeability $K_{eff}$ is Lipschitz-continuous in time and space
\cite{MiWheWi14}.
\end{remark}

\subsection{Spatial discretization of the incremental displacement-phase-field system}

In this section, we formulate a quasi-monolithic Euler-Lagrange formulation 
for $\bU$ and $\Phi$ (approximating $\bu$ and $\phi$), respectively. 
We consider a time-discretized system in which 
time enters through the irreversibility condition. The spatial 
discretized solution variables are 
$\bU \in \mathcal{C}^0([0,T];\pmb{\mathbb{V}}_0(\mathcal{T}))$ and 
$\Phi\in \mathcal{C}^0([0,T];\mathbb{Z}(\mathcal{T}))$, 
where
\begin{align}
&\pmb{\mathbb{V}}_0(\mathcal{T}) := 
\{ W \in C^0(\bar{\Lambda};\mathbb{R}^d) \ | \ W= \boldsymbol{0} \ \text{on } \partial \Lambda, W|_{\mathcal{K}} \in \pmb{\mathbb{Q}}^1(\mathcal{K}), \forall \mathcal{K} \in \mathcal{T} \} , \\
&\mathbb{Z}(\mathcal{T}) := 
\{ Z \in C^0(\bar{\Lambda};\mathbb{R}) | 
\ Z^{n+1}\leq Z^n \leq 1, Z|_{\mathcal{K}} \in \mathbb{Q}^1(\mathcal{K}), \forall \mathcal{K} \in \mathcal{T} \}.
\end{align}
Moreover, we extrapolate $\Phi$ (denoted by $E(\Phi)$) in the first terms
(i.e., the displacement equation) in Formulation \ref{form_2} in order to 
avoid an indefinite Hessian matrix:

\[
E(\Phi) = \Phi^{n-2} + \frac{(t - t^{n-1}-t^{n-2})}{(t-t^{n-1}) -
  (t-t^{n-1}-t^{n-2})} (\Phi^{n-1} - \Phi^{n-2}).
\]
This heuristic procedure has been shown to 
be an efficient and robust method as discussed in \cite{HeWheWi15}. 

In the following, we denote by $\bU^n$ and $\Phi^n$ the approximation of $\bU(t^n)$ and $\Phi(t^n)$ respectively.
\begin{form}
\label{form_2}
Let us assume that $P^{n+1}$ is {a} given approximated pressure at the time $t^{n+1}$.
Given the initial conditions
$\bU^0:=\bU(0)$ and $\Phi^0:=\Phi(0)$ we seek
$\{\bU^{n+1}, \Phi^{n+1} \} \in \pmb{\mathbb{V}}_0(\mathcal{T}) \times \mathbb{Z}(\mathcal{T})$ such that 
\begin{align}
&A_{DPFF}(\bU^{n+1}, \Phi^{n+1})(\bw, \psi-\Phi^{n+1})
\geq 0,  
\ \ \forall \{ \bw, \psi \} \in \pmb{\mathbb{V}}_0(\mathcal{T}) \times \mathbb{Z}(\mathcal{T}),x
\end{align}
with
\begin{align}
&A_{DPFF}(\bU^{n+1}, \Phi^{n+1})(\bw, \psi-\Phi^{n+1})\\ 
&\quad =\int_{\Lambda} (1-k) ( E(\Phi^{n+1})^2 + k) \sigma^+(\bU^{n+1}) : e(\bw) \ d\bx
+\int_{\Lambda} \sigma^-(\bU^{n+1}) : e(\bw) \ d\bx \\
&\qquad - \int_{\Lambda} (\alpha-1) E(\Phi^{n+1})^2 {P^{n+1}} \nabla \cdot \bw  \ d\bx
+ \int_{\Lambda} E(\Phi^{n+1})^2 \nabla {P^{n+1}} \cdot \bw    \ d\bx  \\ 
&\qquad + (1-k)\int_{\Lambda} \Phi^{n+1} \sigma^+(\bU^{n+1}): e(\bU^{n+1}) \cdot (\psi
- \Phi^{n+1})\ d\bx \\
&\qquad - 2(\alpha-1) \int_{\Lambda}\Phi^{n+1} {P^{n+1}} \nabla \cdot \bU^{n+1}
\cdot (\psi - \Phi^{n+1}) \ d\bx
{  + \int_{\Lambda} 2\Phi^{n+1} \nabla {P^{n+1}} \cdot \bU^{n+1} \cdot (\psi  - \Phi^{n+1})   \ d\bx   }
\\
&\qquad  -  G_c \int_{\Lambda}  \dfrac{1}{\varepsilon} ( 1-\Phi^{n+1}) \cdot (\psi- \Phi^{n+1}) \ d\bx
 + G_c \int_{\Lambda} {\varepsilon} \nabla \Phi^{n+1}  \cdot\nabla (\psi- \Phi^{n+1}) \ d\bx.
\end{align}
The solution of this nonlinear variational inequality is briefly 
explained in Section \ref{sec_algo} with all details presented in \cite{HeWheWi15}.
\end{form}

\subsection{Variational formulations of the level-set and the width problems}

The spatially discretized solution variables for 
the level-set and the width are denoted by $\Phi_{\textsf{LS}}(\bx,t)$ and $W(\bx,t)$, respectively. 
Those functions are approximated by using continuous piecewise polynomials
given in their respective finite element spaces,

$$
{\mathbb{V}}_{\textsf{LS}}(\mathcal{T}) := 
\{ \Psi \in C^0(\bar{\Lambda};\mathbb{R}) \ | \ \Psi = {0} \ \text{on }
\Gamma_F, \Psi|_{\mathcal{K}} \in {\mathbb{Q}}^1(\mathcal{K}), \forall \mathcal{K} \in \mathcal{T} \}
$$
for the level-set and 
$$
{\mathbb{V}}_{w}(\mathcal{T}) := 
\{ \Psi \in C^0(\bar{\Lambda};\mathbb{R}) \ | \ \Psi= {0} \ \text{on }
\partial\Lambda, \Psi|_{\mathcal{K}} \in {\mathbb{Q}}^1(\mathcal{K}), \forall \mathcal{K} \in \mathcal{T} \}
$$
for the width.
Assuming that the displacement field $\bU^n$
and the phase field $\Phi^n$ 
are given at time $t^n$,
the Galerkin approximation of the system in Formulation \ref{form:level_set} is formulated as follows:
\begin{form}
\label{form_level_set_discrete}
Find $\Phi_{\textsf{LS}} \in C^0([0,T];{\mathbb{V}}_{\textsf{LS}}(\mathcal{T}))$ such that
\[
A_{\textsf{LS}}(\Phi_{\textsf{LS}},\psi) 
= F_{\textsf{LS}}(\psi) \quad
\forall {\psi\in {\mathbb{V}}_{\textsf{LS}}}(\mathcal{T}),
\]
where
\begin{align*}
A_{\textsf{LS}}(\Phi_{\textsf{LS}},\psi) &:= (\nabla\Phi_{\textsf{LS}},\nabla\psi)
 +  \theta\int_{\Gamma_F^n}  \Phi_{\textsf{LS}}\cdot\psi  ds
,\\
{F_{\textsf{LS}}(\psi)} &:= { \left(\chi^n f_1 + (1-\chi^n) f_2,\psi \right)},
\end{align*}
{ and $\theta \approx 10^3$ is a {sufficiently large} penalty parameter,
  which 
plays a similar role as in discontinuous Galerkin methods (e.g.,
\cite{PieErn12, Riviere2008}).}
\end{form}
Here 
$\Gamma_F^n := \{ \bx \in \Lambda \ | \ |\Phi^n(\bx) - 0.1| \leq \tilde{\epsilon}  \}$, with a small positive constant $\tilde{\epsilon}$   
and we prescribe this surface with the help of so-called material ids for each cell.
These are set to $0$ for $\Phi^n < C_{LS}$ and $1$ otherwise to identify the interface ($\Gamma_F^n$) on the  discrete level.
It follows that $\chi^n = 0$ for $\Phi^n < C_{LS}$ and $\chi^n = 1$
otherwise. The Galerkin approximation of the width system in Formulation \ref{form:width} is given by
\begin{form}
\label{form_width_discrete}
Find $W\in C^0([0,T];{\mathbb{V}}_w(\mathcal{T}))$ such that
\[
A_{W}(W,\psi) = F_W(\psi) \quad
\forall \psi\in {\mathbb{V}}_{w}(\mathcal{T})
\]
where  
\begin{align*}
A_{W}({W},\psi) &= 
(\nabla {W},\nabla\psi)  + \theta \int_{\Gamma_F^n}
{W} \psi \, ds,\\
F_W(\psi) &= \theta \int_{\Gamma_F^n}
W_D^n \cdot \psi \, ds.
\end{align*}
Here $W_D^n := -2 \bU^n \cdot \frac{\nabla\Phi^n_{\textsf{LS}}}{\|\nabla\Phi^n_{\textsf{LS}}\|}$ is the width on the fracture boundary $\Gamma_F^n$.
\end{form}

\section{Solution algorithms for fluid-filled phase-field fractures}
\label{sec_fs_splitting}
In this section, we formulate iterative coupling of the 
two physical subproblems; namely pressure and displacement/phase-field.
However before each pressure solve we also solve first 
the level-set problem and the width problem in order to 
provide the fracture permeability.

\subsection{Solution algorithms and solver details}
\label{sec_algo}
In Algorithm \ref{alg:solve_algo_fixedstress},
 we outline the entire scheme for all solution variables 
$\{\Phi_{\textsf{LS}}^l,W^l,P^l,\bU^l,\Phi^l\}$ at each fixed-stress iteration 
step. 

\begin{algorithm}
\caption{Iterative coupling for fluid-filled phase-field fractures including
  level-set and width computation}
\label{alg:solve_algo_fixedstress}
\begin{algorithmic}
\STATE At the time $t^{n}$,
\REPEAT  
\STATE For $l=0,1,2,\ldots$: 
\STATE -  Solve the (linear) level-set Formulation
\ref{form_level_set_discrete} for $\Phi_{\textsf{LS}}^l$
\STATE -  Solve the (linear) width Formulation \ref{form_width_discrete} for $W^l$
\STATE -  Solve the (linear) pressure diffraction Formulation \ref{form_1} for $P^l$
\STATE -  Solve the (nonlinear) fully-coupled displacement/phase-field
Formulation \ref{form_2} for $(\bU^l, \Phi^l)$
\UNTIL
the stopping criterion for fixed-stress split is satisfied:
\begin{equation*}
\max\{ \|P^l - P^{l-1}\| , \| \bU^l - \bU^{l-1} \| , \| \Phi^l - \Phi^{l-1} \|\} \leq \operatorname{TOL_{FS}},
\quad\operatorname{TOL_{FS}}>0
\end{equation*}
\STATE Set: $(P^n, \bU^n,\Phi^n):= (P^l, \bU^l, \Phi^l)$. The other two 
variables $\Phi_{\textsf{LS}}^n$ and $W^n$ are obtained from $\Phi^n$ and $\bU^n$.
\STATE Increment the time $n\to n+1$.
\end{algorithmic}
\end{algorithm}
\begin{remark}
{In Algorithm \ref{alg:solve_algo_fixedstress}, the explicit 
solution of $\Phi_{\textsf{LS}}^l$ is avoided when working with 
Formulation \ref{form_level_set_by_PFF}.}
\end{remark}

The nonlinear quasi-monolithic displacement/phase-field 
system presented in Formulation \ref{form_2} is solved with Newton's method and line search
algorithms. The constraint minimization problem is treated with 
a semi-smooth Newton method (i.e., a primal-dual active set method). 
Both methods are combined in one single loop leading to a robust and efficient 
iteration scheme that is outlined in  \cite{HeWheWi15}. 
Within Newton's loop we solve the linear equation systems 
with GMRES solvers with diagonal block-preconditioning from Trilinos 
\cite{Trilinos-Overview}. 
Algorithm \ref{alg:solve_algo_fixedstress} presents the overall fixed-stress 
phase field approach for fluid filled fractures in which the
geomechanics-phase-field system is coupled to the pressure diffraction problem. 
We employ local mesh adaptivity in order to keep the computational cost at a
reasonable level. Here, we specifically use a technique developed 
in \cite{HeWheWi15}; namely a predictor-corrector scheme  
that chooses an initial $\eps > h$ at the beginning of the computation.
Since this is a model parameter, we do not want to change the 
model during the computation and keep therefore $\ep$ fixed. 
However, the crack propagates and in coarse mesh regions $\eps > h$ may
 be violated. Then, we take the first step as a predictor step,
then refine the mesh such that $\eps > h$ holds again and finally recompute 
the solution. In \cite{HeWheWi15} (two dimensional) and \cite{WiLeeWhe15,LeeWheWi16} (three dimensional) 
it has been shown that this procedure is efficient and robust.
The pressure diffraction problem is solved with a direct solver for simplicity
but could also be treated with an iterative solver.
The linear-elliptic level-set and width problems are solved 
with a parallel CG solver and SSOR preconditioning where the 
relaxation parameter is chosen as $1.2$.

\subsection{A fixed-stress algorithm for fluid-filled phase-field fractures}
In this section, we now focus on the specifics of the fixed-stress iteration 
between flow and geomechanics/fracture.
Let $\Phi^l$ and $W^l$ at time $t^{n+1}$ be given.
For each time $t^{n+1}$ we iterate for $l=0,1,2,3,\ldots$:
\paragraph{i) Fixed-stress: pressure solve}
Let $\bU^l$ and $\Phi^l$ be given. Find $P^{l+1} \in \mathbb{W}$ such that
$$
[A_P(P^{n+1})(\omega)]^{l+1} = [A_{PR}(P^{n+1})(\omega)]^{l+1} + [A_{PF}(P^{n+1})(\omega)]^{l+1} = 0 \quad\forall\omega\in\mathbb{W}(\mathcal{T}),
$$
where
\begin{multline}
[A_{PR}(P^{n+1})(\omega)]^{l+1} := 
\chi_R(\Phi^{l}) \Biggl( \int_{\Lambda}\rho_R^0 \Big(\dfrac{1}{M} + \frac{3\alpha^2}{3\lambda + 2\mu} \Big) \Big( \dfrac{P^{l+1} - P^n}{\delta t}  \Big) \cdot \omega \ d\bx
+ \int_{\Lambda}\dfrac{K_R \rho_R^0}{\eta_R} (\nabla P^{l+1} - \rho_R^0 \bg) \nabla \omega \ d\bx  \\
+\int_{\Lambda} \alpha \nabla \cdot  \Big( \dfrac{\bU^{l} -\bU^{n}}{\delta t}  \Big) \cdot \omega \ d\bx 
-\int_{\Lambda} \Big(\frac{3\alpha^2}{3\lambda + 2\mu} \Big) \Big( \dfrac{P^{l} - P^n}{\delta t}  \Big) \omega \ d\bx
-\int_{\Lambda} q_R  \omega  \ d\bx \Biggr),
\quad \forall \omega \in \mathbb{W}(\mathcal{T}),
\end{multline}
\begin{multline}
[A_{PF}(P^{n+1})(\omega)]^{l+1} :=
\chi_F(\Phi^{l})  \Biggl( \int_{\Lambda} \rho_F^0 c_F \Big( \dfrac{P^{l+1} - P^n}{\delta t}  \Big)  \omega \ d\bx
+ \int_{\Lambda} \dfrac{K_F \rho_F^0}{\eta_F} (\nabla P^{l+1} - \rho_F^0 \bg)  \nabla \omega \ d\bx  \\
- \int_{\Lambda} q_F \omega \ d\bx \Biggr), 
\quad \forall \omega \in \mathbb{W}(\mathcal{T}).
\end{multline}

\paragraph{ii) Fixed-stress: displacement/phase-field solve}
Take the just computed $P^{l+1}$ and solve for the displacements $\bU^{l+1}\in
\pmb{\mathbb{V}}_0(\mathcal{T})$ and the phase field $\Phi^{l+1} \in
\mathbb{Z}(\mathcal{T})$ such that:
\begin{equation}
A_{DPFF}(\bU^{l+1},\Phi^{l+1})(\bw,\psi) \geq 0 \quad \forall \{ \bw, \psi \} \in \pmb{\mathbb{V}}_0(\mathcal{T}) \times \mathbb{Z}(\mathcal{T}),
\end{equation}
where
\begin{align*}
A_{DPFF}(\bU^{l+1}, \Phi^{l+1})(\bw, \psi) &=
\int_{\Lambda}( (1-k) ( E(\Phi^{l+1})^2 + k) \sigma^+(\bU^{l+1}) : e(\bw) \ d\bx
+\int_{\Lambda} \sigma^-(\bU^{l+1}) : e(\bw) \ d\bx \\
&\quad - \int_{\Lambda} (\alpha-1) E(\Phi^{l+1})^2 P^{l+1}  
\nabla \cdot \bw  \ d\bx
+ \int_{\Lambda} E(\Phi^{l+1})^2 \nabla P^{l+1} \cdot \bw    \ d\bx\\   
&\quad + (1-k)\int_{\Lambda} \Phi^{l+1} \sigma^+(\bU^{l+1}): e(\bU^{l+1}) \cdot \psi \ d\bx\\ 
&\quad - 2(\alpha-1) \int_{\Lambda}\Phi^{l+1} P^{l+1} \nabla \cdot \bU^{l+1} \cdot\psi \ d\bx
{  + \int_{\Lambda} 2\Phi^{l+1} \nabla P^{l+1} \cdot \bU^{l+1} \cdot\psi     \ d\bx   }
\\
&\quad -  G_c \int_{\Lambda}  \dfrac{1}{\varepsilon} ( 1-\Phi^{l+1}) \cdot \psi \ d\bx
 + G_c \int_{\Lambda} {\varepsilon} \nabla \Phi^{l+1}  \cdot\nabla \psi \ d\bx.
\end{align*}
\paragraph{iii) Fixed-stress: stopping criterion} The iteration is completed if 
\[
\max\{ \| \bU^{l+1} - \bU^{l} \|_{L^2(\Lambda)}, \|P^{l+1} - P^{l}\|_{L^2(\Lambda)}, \|\Phi^{l+1} - \Phi^{l}\|_{L^2(\Lambda)} \} < TOL_{FS}.
\]
then we set
$$
P^{n+1} = P^{l+1}, \ \Phi^{n+1} = \Phi^{l+1}, \ \bU^{n+1} = \bU^{l+1}.
$$
The specific tolerances $TOL_{FS}$ will be specified in Section \ref{sec_tests}.

\section{Numerical tests}
\label{sec_tests}
In this final section, we present five different examples with 
increasing complexity. Our main focus is on fixed-stress 
iteration numbers and refinement studies in order 
to investigate the iterative solution approach.
{ Alongside we discuss differences between $\alpha=0$
and $\alpha = 1$, negative pressure at fracture tips, 
the pressure drop for propagating fractures,
and interaction of multiple fractures in heterogeneous media.}
The examples are computed with the finite element 
package deal.II \cite{BangerthHartmannKanschat2007,dealII84} 
and are based on the 
programming codes developed in \cite{HeWheWi15,WiLeeWhe15,LeeWheWi16}
by using an MPI-parallel framework. 

\paragraph{Boundary and initial conditions}
For all following examples, the initial crack is given with the help of the phase-field function $\varphi$. We set at $t=0$:
\begin{equation}
\varphi = 0 \quad\text{in } \Omega_F, \quad\text{and } \quad
\varphi = 1 \quad\text{in } \Lambda\setminus\Omega_F
\end{equation}
for each defined $\Omega_F$.
As boundary conditions, we set the
displacements to zero  on $\partial\Omega$ and 
traction-free conditions for the phase-field variable.
The boundary and interface conditions for the pressure,
level-set and width computation have been explained 
in their respective sections before.
In addition, we recall that 
the diameter of an element $\mathcal{K}\in \mathcal{T}_h$ is denoted by $h$,
and $h_{\min}$ for the minimum diameter, 
and $h_{\max}$ for the maximum diameter during adaptive mesh refinement.

\subsection{Example 1: Extension of Sneddon's test to a fluid-filled fracture
  in a porous medium}
Sneddon's test {\cite{Snedd46,SneddLow69}} 
is an important example for pressurized fractures in which 
a given pressure causes the fracture to open. The pressure is 
however too low to propagate the fracture in its length. 
In this first example, we extend this test to a fluid-filled setting
in which fluid is injected into the middle of the fracture.
The flow rate injection $q_F$ is chosen as such that Sneddon's pressure (for 
example $p\approx 10^{-3}$ as used in \cite{BourChuYo12,WheWiWo14}) 
is approximately recovered and to study whether the resulting crack
opening displacement is of the same order as in the existing 
literature. Here we then carry out convergence studies with 
respect to the mesh size parameter $h$ while keeping the 
model parameter $\ep$ fixed.

\paragraph{Configuration}
We deal with the following geometric data:
$\Omega = (\SI{0}{\metre},\SI{4}{\metre})^2$ and a (prescribed) initial crack  
with half length $l_0 = \SI{0.2}{\metre}$ 
on $\Omega_F = (1.8,2.2)\times(2-h_{\max}, 2+h_{\max}) \subset \Omega $. 
The initial mesh is $5$ times uniformly refined, 
and then $3,4$ and $5$ times locally, sufficiently large 
around the fracture region. 
This leads to 
$3580, 10312$ and $36052$  initial mesh cells, with 
$h_{min}=\SI{0.022}{\metre}, \SI{0.011}{\metre}$ 
and $\SI{0.0055}{\metre}$, respectively.
On the finest mesh we have 
73018 degrees of freedom (DoFs) for the solid, and $36509$ DoFs for each the phase-field,
the pressure, the level-set and the width, respectively.

\paragraph{Parameters}
The critical energy release rate is chosen as
$G_c = \SI{1}{\newton\per\metre}$. The mechanical parameters are
Young's modulus and Poisson's ratio $E = \SI{1}{\pascal}$ and $\nu_s = 0.2$.
The relationship to the Lam\'e coefficients $\mu_s$ and $\lambda_s$ is
given by:
\begin{align*}
\mu_s &= \frac{E}{2(1+\nu_s)} , \qquad 
\lambda_s = \frac{\nu_s E_s}{(1+\nu_s )(1-2\nu_s)}. 
\end{align*}
{The regularization parameters are chosen as $\eps = 2h_{max} = 0.045$
and $\kappa = \num{e-10}h_{min}$}. 
We perform computations for Biot's coefficient $\alpha = 0$ and $\alpha = 1$.
Furthermore $q_F = \SI{5e-2}{m^3/s}$ for $\alpha = 1$ and $q_F =
\SI{5e-9}{m^3/s}$ for $\alpha = 0$, 
and $M = \SI{1e+12}{Pa}, c_F = \SI{1e-12}{Pa^{-1}}$. The viscosities are chosen as
$\eta_R = \eta_F = \SI{1e-3}{Ns/m^2}$. The reservoir permeability is $K_R =
\SI{1}{d} = \SI{1e-12}{m^2}$  
and the density is $\rho_F^0 = \SI{1}{kg/m^3}$.
Furthermore, $TOL_{FS} = 10^{-3}$.
This test case is computed in a quasi-stationary manner, which 
is due to the crack irreversibility constraint. That is,
we solve $10$ pseudo-time steps with a time step size $k=1s$.

\paragraph{Quantities of interest}
We study the following cases and goal functionals:
\begin{itemize}
\item the maximum pressure evolution over time; 
\item the crack opening displacement (or aperture):
\begin{equation}
\label{Jump2}
COD(x_0) := [\bu(x_0,y)\cdot\bn] 
= \int_{0}^{4}  \bu(x_0, y)\cdot \bn \, dy 
= \int_{0}^{4}  \bu(x_0, y)\cdot \nabla\varphi(x_0,y) \, dy,
\end{equation}
where $\varphi$ is our phase-field function and 
$x_0$ the $x$-coordinate of the integration line.
The analytical solution for the crack opening displacement is 
derived by Sneddon and Lowengrub \cite{SneddLow69};
\item the number of fixed-stress iterations;
\item the number of Newton solves for the displacement/phase-field system;
\item the number of GMRES iterations within the Newton solver.
\end{itemize}

\paragraph{Discussion of findings}
We observe $4$ fixed-stress iterations in the very first time step 
for reaching the tolerance $TOL_{FS} = 10^{-3}$. In the subsequent 
time steps, we immediately satisfy the tolerance in one step since 
the problem is quasi-stationary and not much change happens.
To solve the nonlinear displacement/phase-field system, 
$2-4$ Newton steps (that then satisfy both 
criteria: active set convergence and the nonlinear residual tolerance) 
are required in average. Inside each Newton 
step we need in average $10-40$ GMRES iterations. Here we do not 
observe significant differences between $\alpha = 0$ and $\alpha = 1$.

\begin{figure}[H]
\centering
\begin{subfigure}[b]{0.35\textwidth}
{\includegraphics[width=\textwidth]{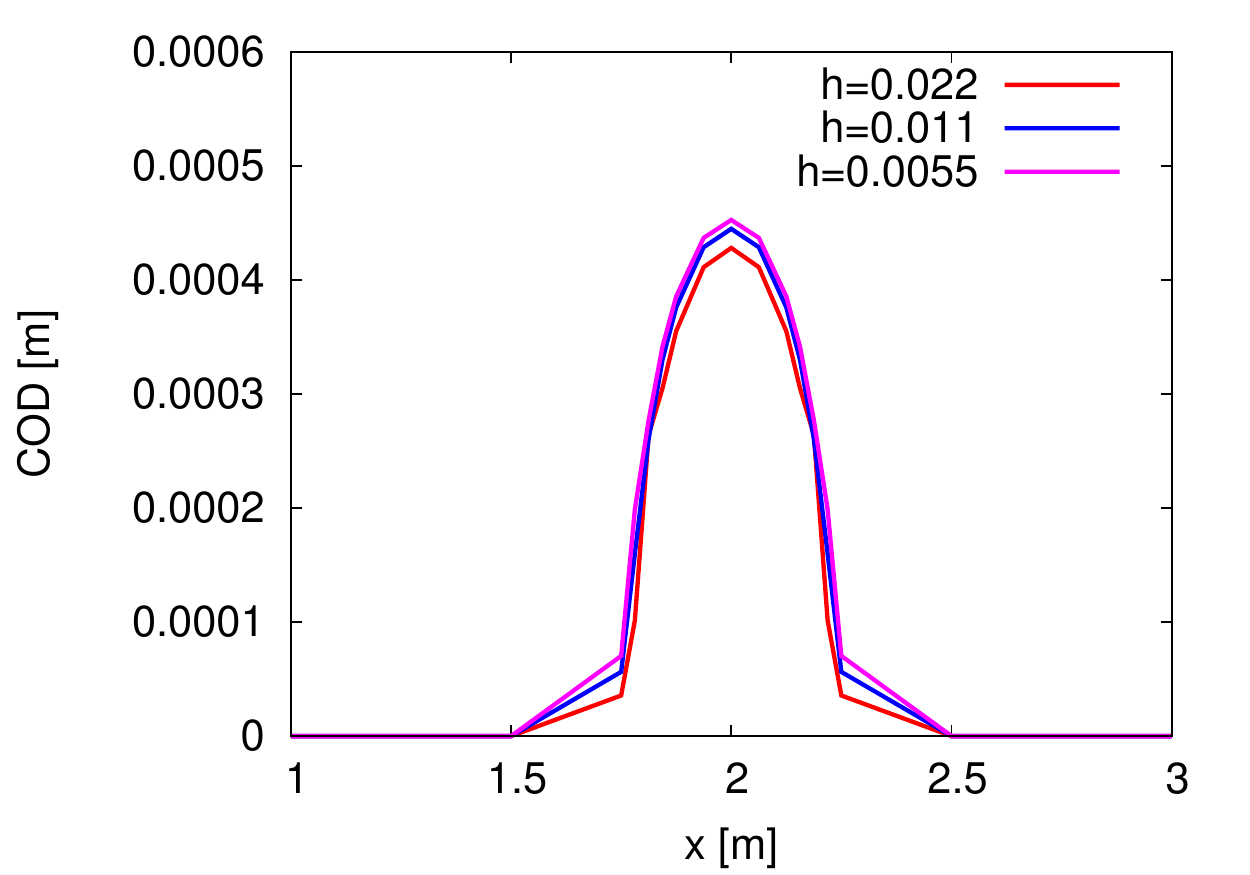}}
\caption{$\alpha=0$}
\end{subfigure}
\begin{subfigure}[b]{0.35\textwidth}
{\includegraphics[width=\textwidth]{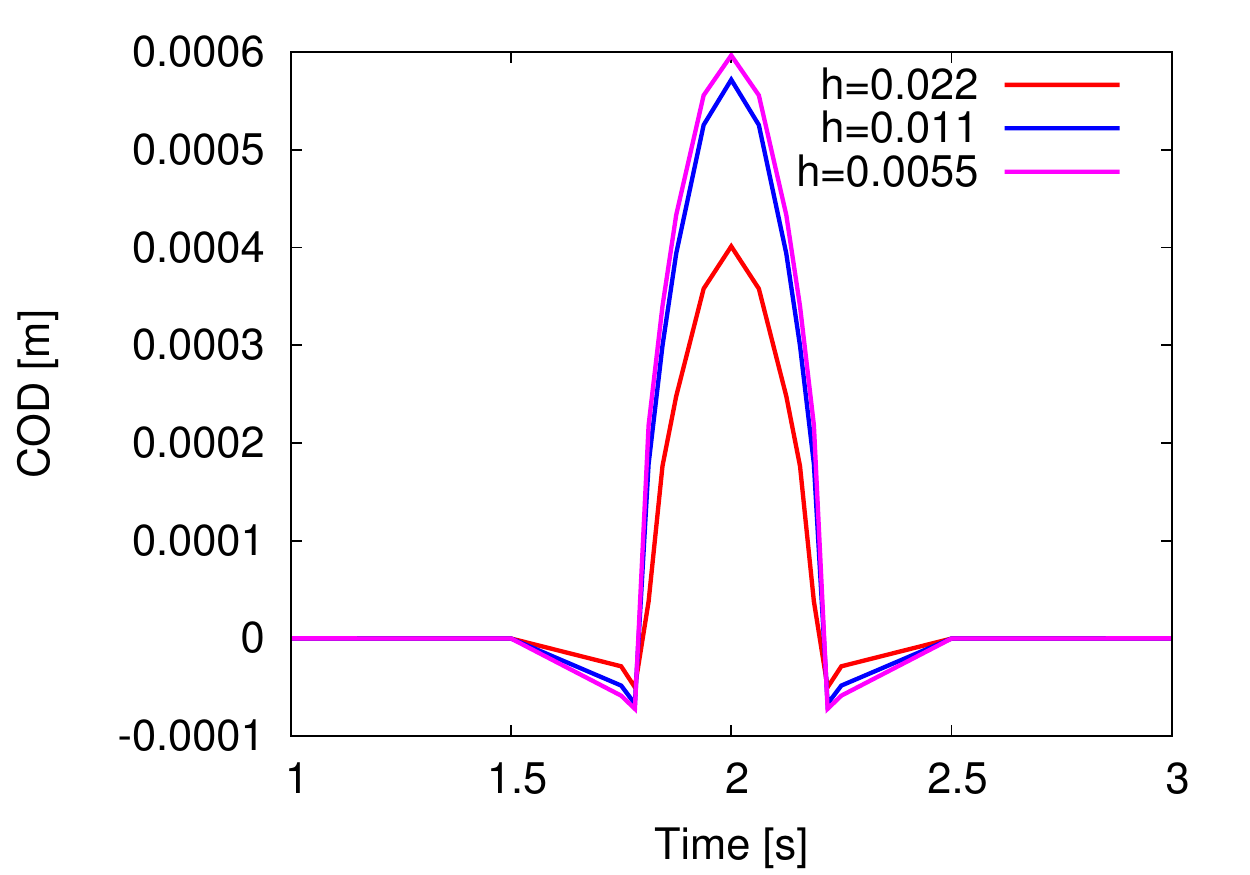}}
\caption{$\alpha=1$}
\end{subfigure}
\begin{subfigure}[b]{0.35\textwidth}
{\includegraphics[width=\textwidth]{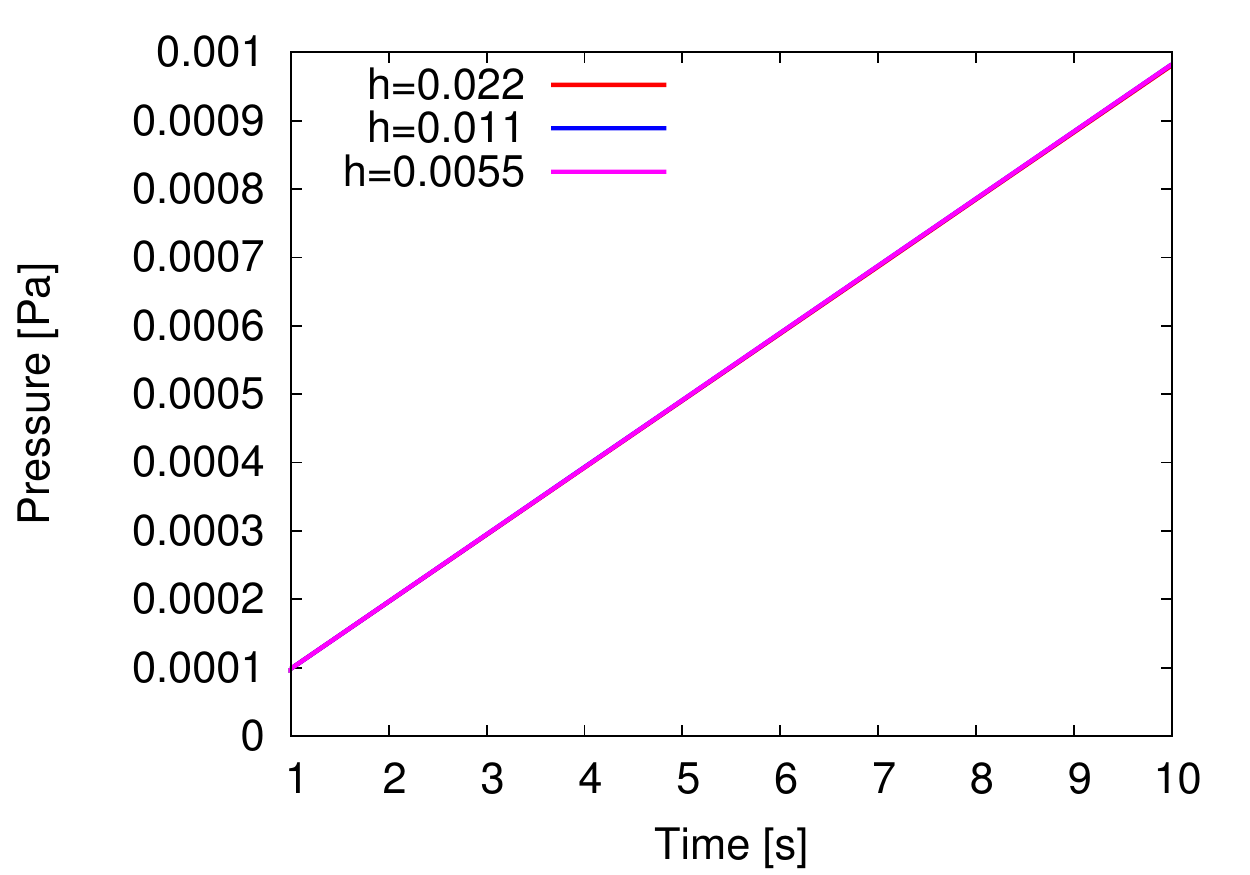}}
\caption{$\alpha=0$}
\end{subfigure}
\begin{subfigure}[b]{0.35\textwidth}
{\includegraphics[width=\textwidth]{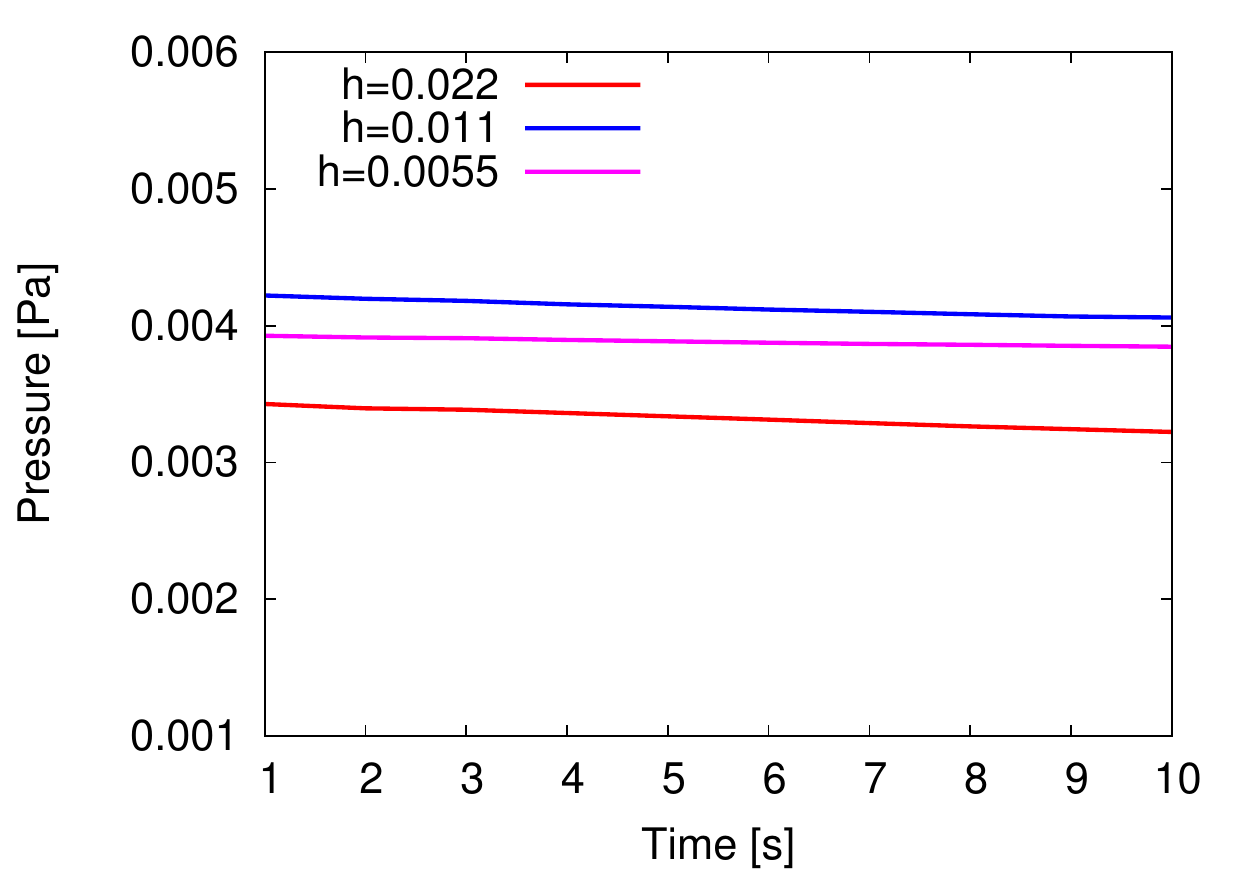}}
\caption{$\alpha=1$}
\end{subfigure}
\caption{Example 1: On top: the crack opening displacement for (a) $\alpha = 0$ and (b)  $\alpha
  = 1$ at $T=10s$. 
On the bottom: (c) - (d) the corresponding maximal pressure evolutions. 
{In the bottom left subfigure at $T=10s$ the desired pressure $p\approx 0.001
Pa$ is reached. This test ($\alpha = 0$) compares (as expected) very much 
to the original Sneddon's test with a given, fixed pressure. 
The test with $\alpha = 1$ differs since now poroelastic effects play a role 
that were not accounted for in Sneddon's original derivation.}}  
\label{ex_1_fig_1}
\end{figure}

\begin{figure}[!h]
\centering
\begin{subfigure}[b]{0.225\textwidth}
{\includegraphics[width=\textwidth]{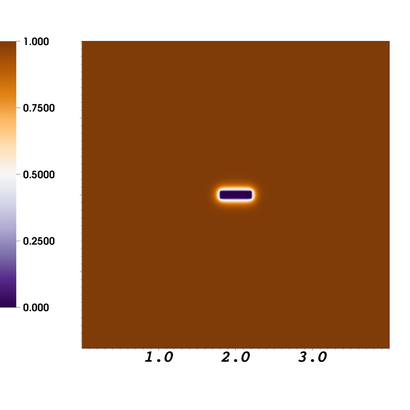}}
\caption{phase-field $\varphi$}
\end{subfigure}
\begin{subfigure}[b]{0.225\textwidth}
{\includegraphics[width=\textwidth]{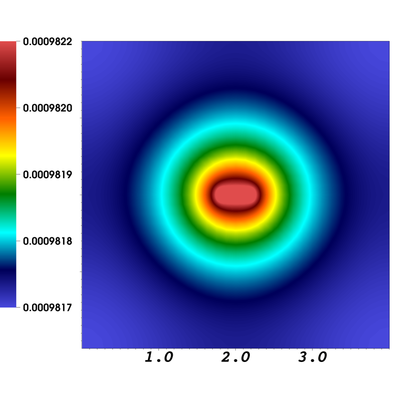}}
\caption{pressure $p$}
\end{subfigure}
\begin{subfigure}[b]{0.225\textwidth}
{\includegraphics[width=\textwidth]{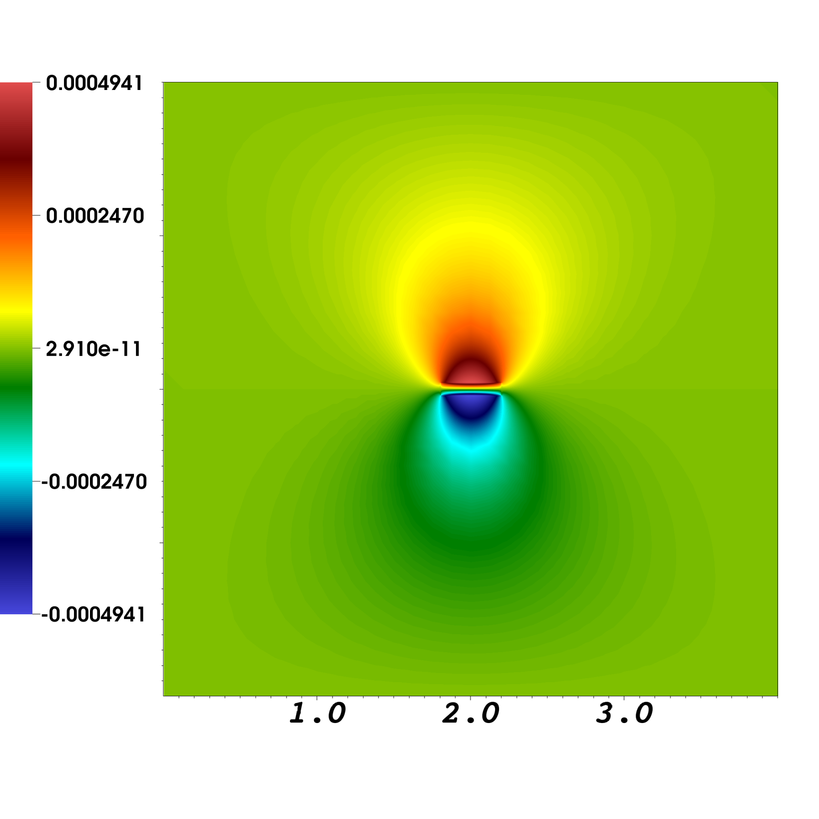}}
\caption{$\bu_y$ displacement}
\end{subfigure}
\begin{subfigure}[b]{0.225\textwidth}
{\includegraphics[width=\textwidth]{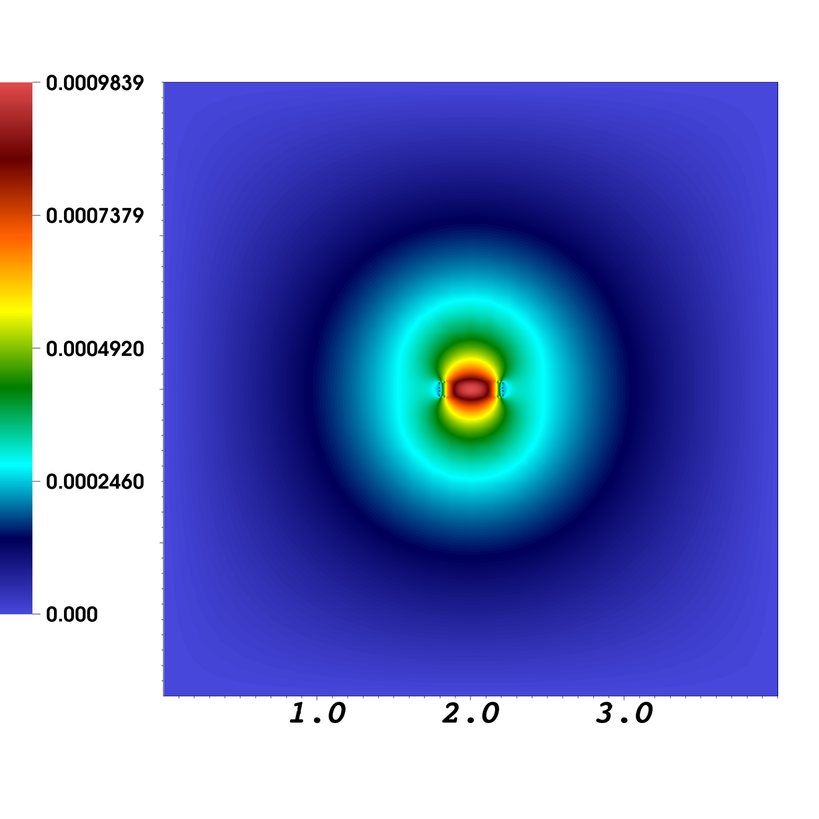}}
\caption{interpolated width $w$}
\end{subfigure}
\caption{Example 1: Case $\alpha = 0$: 
The phase-field $\varphi$, 
the pressure $p$, the $\bu_y$ displacement, and the
interpolated width $w$ are presented.
}
\label{ex_1_fig_2}
\end{figure}

\begin{figure}[!h]
\centering
\begin{subfigure}[b]{0.225\textwidth}
{\includegraphics[width=\textwidth]{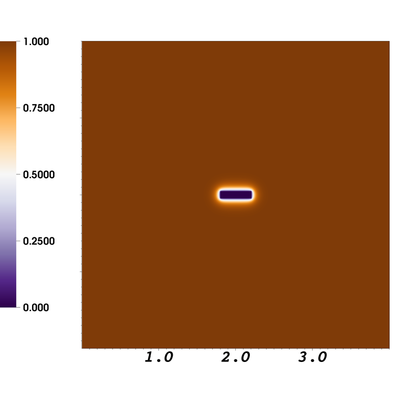}}
\caption{phase-field  $\varphi$}
\end{subfigure}
\begin{subfigure}[b]{0.225\textwidth}
{\includegraphics[width=\textwidth]{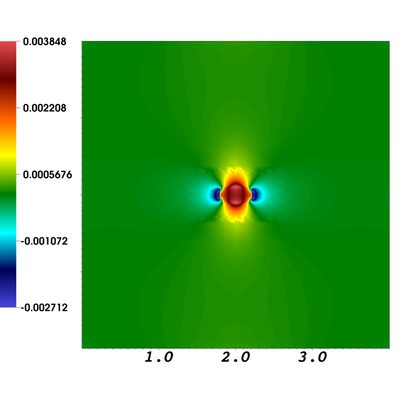}}
\caption{pressure $p$}
\end{subfigure}
\begin{subfigure}[b]{0.225\textwidth}
{\includegraphics[width=\textwidth]{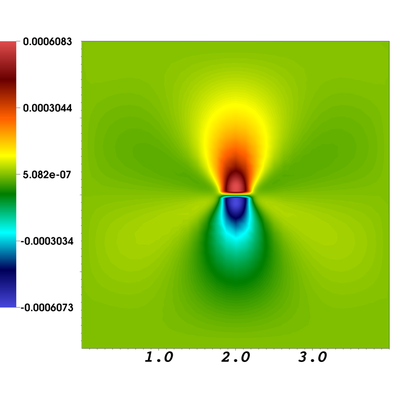}}
\caption{$\bu_y$ displacement}
\end{subfigure}
\begin{subfigure}[b]{0.225\textwidth}
{\includegraphics[width=\textwidth]{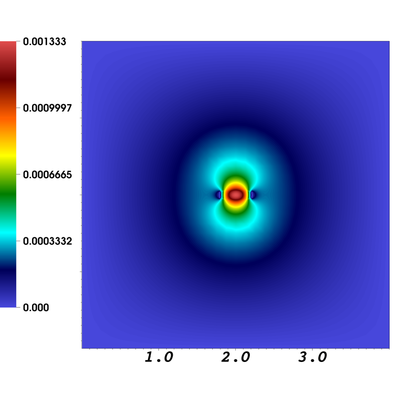}}
\caption{interpolated width $w$}
\end{subfigure}
\caption{Example 1: Case $\alpha = 1$: 
The phase-field $\varphi$, the pressure $p$, the $\bu_y$ displacement, and the
  interpolated width $w$ are presented.}
\label{ex_1_fig_3}
\end{figure}

In Figure \ref{ex_1_fig_1}, we observe mesh convergence 
of the crack opening displacement and the corresponding 
pressure evolutions for both cases $\alpha = 0$ and $\alpha = 1$.
We see that using $\alpha = 0$ the pressure is $p\approx \SI{1e-3}{Pa}$ as used 
by \cite{BourChuYo12,WheWiWo14} and yields a similar crack opening 
displacement. This is a major achievement that the 
fluid-filled model (namely pressure diffraction coupled to
displacement/phase-field) is able to represent the manufactured solution
of the original pressurized test case.
Using $\alpha = 1$ we observe that the pressures are $4$ times higher 
(yielding a slightly higher COD). This is expected since 
in the $\alpha = 1$-case the fracture pressure interacts with 
the reservoir pressure and fluid is released into the porous medium. 
In the Figures \ref{ex_1_fig_2} and \ref{ex_1_fig_3} the different 
solution variables at the end time value $T=10s$ are displayed.
In fact we nicely identify the interpolated width $w$ and also 
the different pressure distributions depending on the different 
$\alpha$ choices.
Moreover, observing the quantitative values for $\bu_y$ (which corresponds 
in this symmetric test to the COD) and the subsequent FE width value, we 
identify excellent agreement.

\subsection{Example 2: A fluid-filled fracture with emphasis on the pressure at the fracture tips.}
\label{sec_ex_2}
In this short section, we only focus on the pressure evolution for a fluid-filled (namely choosing $\alpha=1$) configuration.
We highlight negative pressure values at the fracture tips, which are known as a typical phenomena caused by fluid lagging in the early injection stages and has
been observed by others as well, see for example \cite{DeGara03,Markert2015,SavitskiDetournay2002,SchrefSecSi06,SecSchref12}.

\paragraph{Configuration}
In the domain $\Omega=(-\SI{5}{\metre},\SI{5}{\metre})^2$ the initial penny
shaped fracture is given in the center $(0,0)$ with the longer radius
$r=\SI{1}{\metre}$ on $\Omega_F = (4,6)\times(5-h_{\max}, 5+h_{\max}) \subset \Omega $.
The physical parameters are chosen as $\alpha=1$,
$E=\SI{1}{\pascal}, \nu = 0.2$, $G_c = \SI{1.}{\newton\per\metre}$,
$\eta_R = \eta_F = \SI{1e-3}{Ns/m^2}$, $\rho=\SI{1}{kg/m^3}$, 
$K_R = \SI{1e-16}{m^2}$  
and
$q_F=\SI{7.5e-2}{m^3/s}$.  
Here the numerical parameters are given as $h_{\min}=\SI{0.028}{\metre}$, $\delta t=\SI{0.5}{s}$, $T=\SI{30}{s}$, and $\varepsilon = 2 h_{\min}$.

\begin{figure}[!h]
\centering
\includegraphics[scale=0.32]{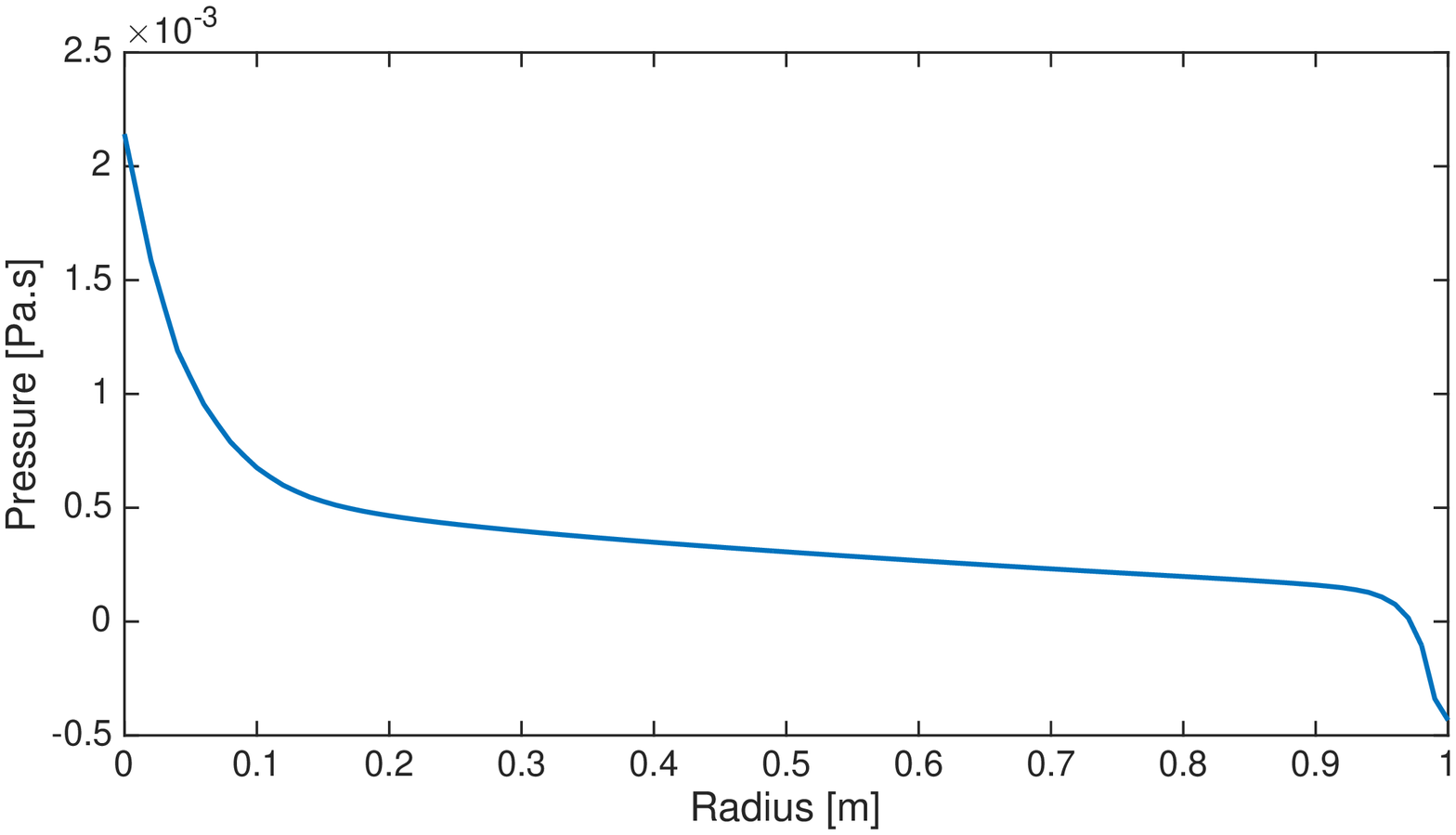}
\includegraphics[scale=0.2]{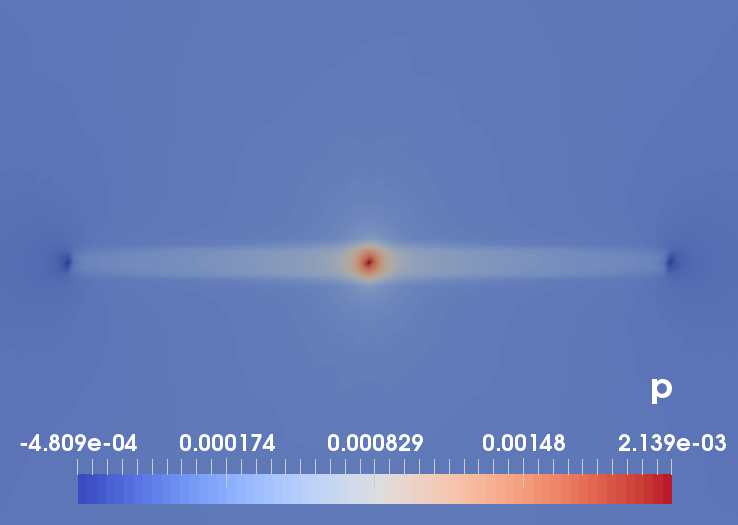}
\caption{Example 2. 
Illustration of the pressure values a) inside (plotted from the center to the end-tip) and b) around the fracture. We observe negative pressures at the fracture tip.}
\label{fig:2}
\end{figure}

\paragraph{Discussion of our findings}
At the early stage of the fluid injection, Figure \ref{fig:2} illustrates the pressure values in the fracture. 
We notice that our findings show a qualitative similarity 
with the plots presented in \cite{SavitskiDetournay2002,Lecampion:2013fr}.
However, due to $\alpha=1$,
the full coupling of a fluid-filled fracture with the surrounding porous medium, and the use of another well model, 
we cannot expect a quantitative agreement
with \cite{SavitskiDetournay2002,Lecampion:2013fr}.
However, the important result is that (as predicted in the above mentioned papers)
that we observe negative pressures around the fracture tips. This 
phenomena may arise in the case of injections if the speed at which the crack
tip advances is sufficiently high such that the fluid inside the fracture 
cannot flow fast enough to fill the created space. In particular, 
at the beginning of an injection, the fracture is not yet completely 
filled with fluid and thus at the tips fluid enters from the porous 
medium into the fracture causing the predicted negative pressure values.

\subsection{Example 3: A propagating fluid-filled fracture in a porous medium}
In this third example, we consider a single propagating fracture. 
The main purpose is to study fixed-stress iterations for this 
nonstationary case.

\paragraph{Configuration}
We deal with the following geometric data:
$\Omega = (\SI{0}{\metre},\SI{4}{\metre})^2$ and a (prescribed) initial crack  
with half length $l_0 = \SI{0.2}{\metre}$ 
on $\Omega_F = (1.8,2.2)\times(2-h_{\max}, 2+h_{\max}) \subset \Omega $.
The initial mesh is $5$ times uniformly refined, 
and then $2,3,4$ and $5$ times locally using predictor-corrector 
mesh refinement, sufficiently large 
around the fracture region. The number of mesh cells will grow 
during the computation due to predictor-corrector 
mesh refinement.

\paragraph{Parameters}
The fracture toughness is chosen as
$G_c = \SI{1}{\newton\per\metre}$. The mechanical parameters are
Young's modulus and Poisson's ratio $E = \SI{e+8}{\pascal}$ and $\nu_s = 0.2$.
The regularization parameters are chosen as $\eps = 2h$
and $\kappa = \num{e-10}h$.
We perform computations for Biot's coefficient $\alpha = 1$ only.
Furthermore the injection rate is chosen as $q_F = \SI{2}{m^3/s}$;
and $M = \SI{1e+8}{Pa}, c_F = \SI{1e-8}{Pa^{-1}}$. The viscosities are chosen as
$\nu_R = \nu_F = \SI{1e-3}{Ns/m^2}$. The reservoir permeability is $K_R = \SI{1}{d}$
and the density is $\rho_F^0 = \SI{1}{kg/m^3}$.
Furthermore, $TOL_{FS} = 10^{-3}$.
The total time is $T=\SI{0.6}{s}$. We also perform time convergence 
studies and use as time steps $\delta t=0.01s,0.005s,0.0025s,0.00125s$.
Thus, $60,120,240$ and $480$ time steps are computed, respectively.

\begin{figure}[H]
\centering
{\includegraphics[width=7cm]{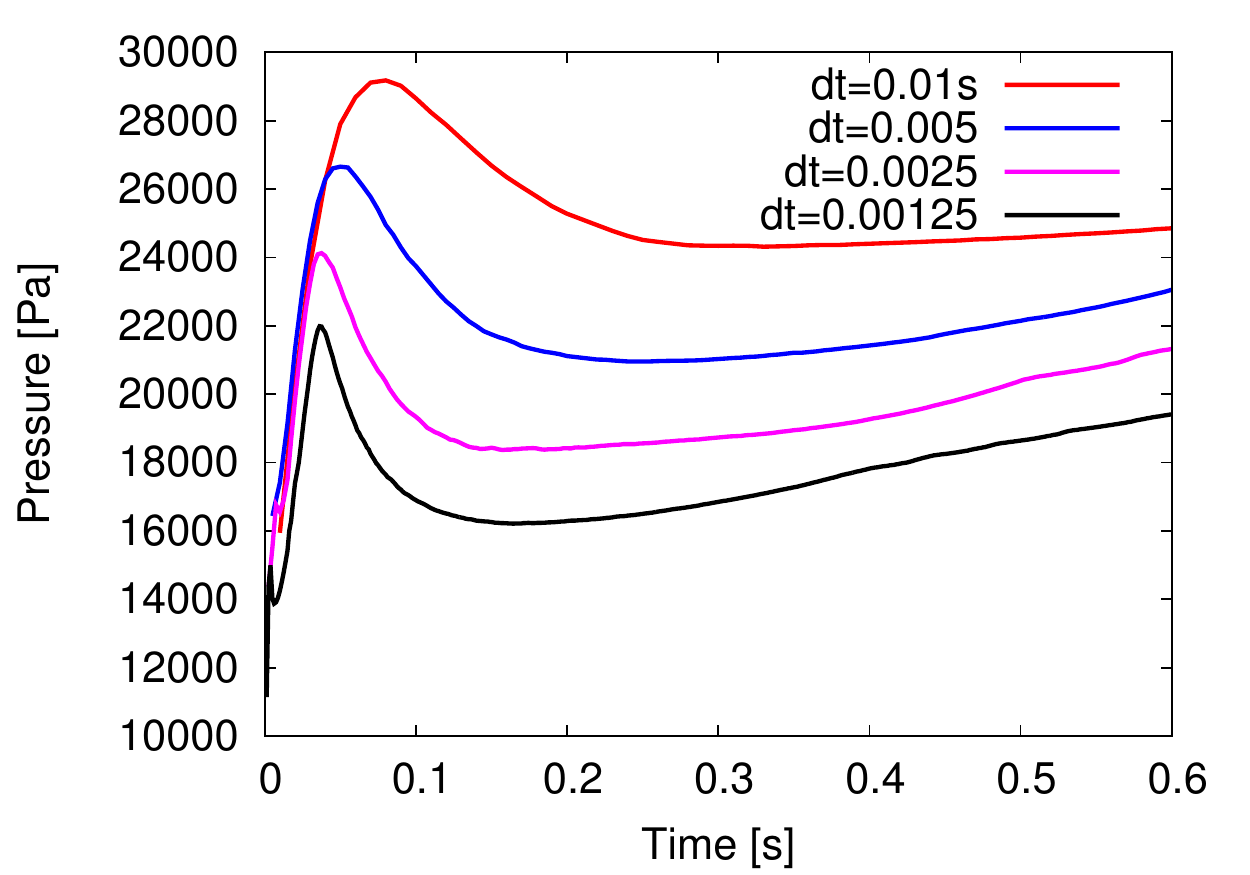}}
{\includegraphics[width=7cm]{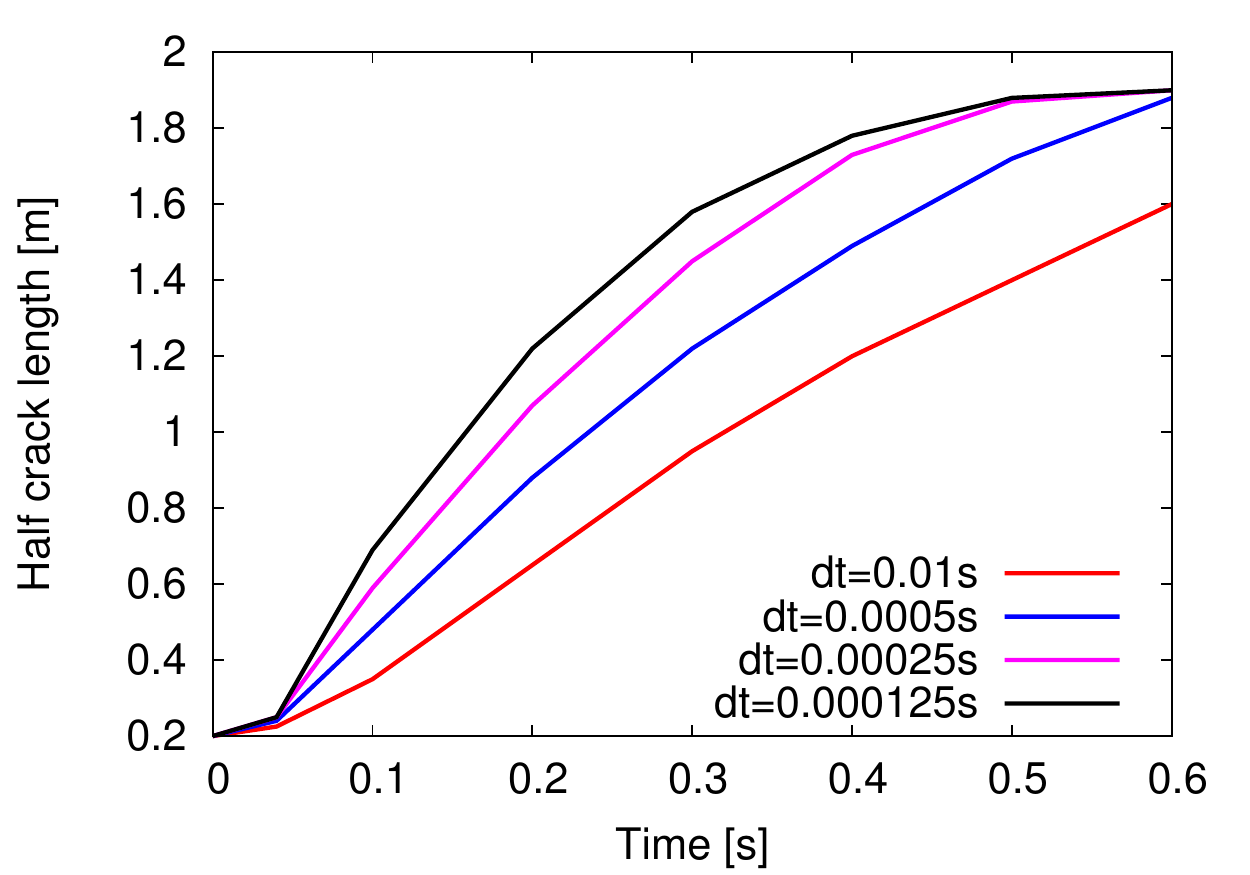}}\\
{\includegraphics[width=7cm]{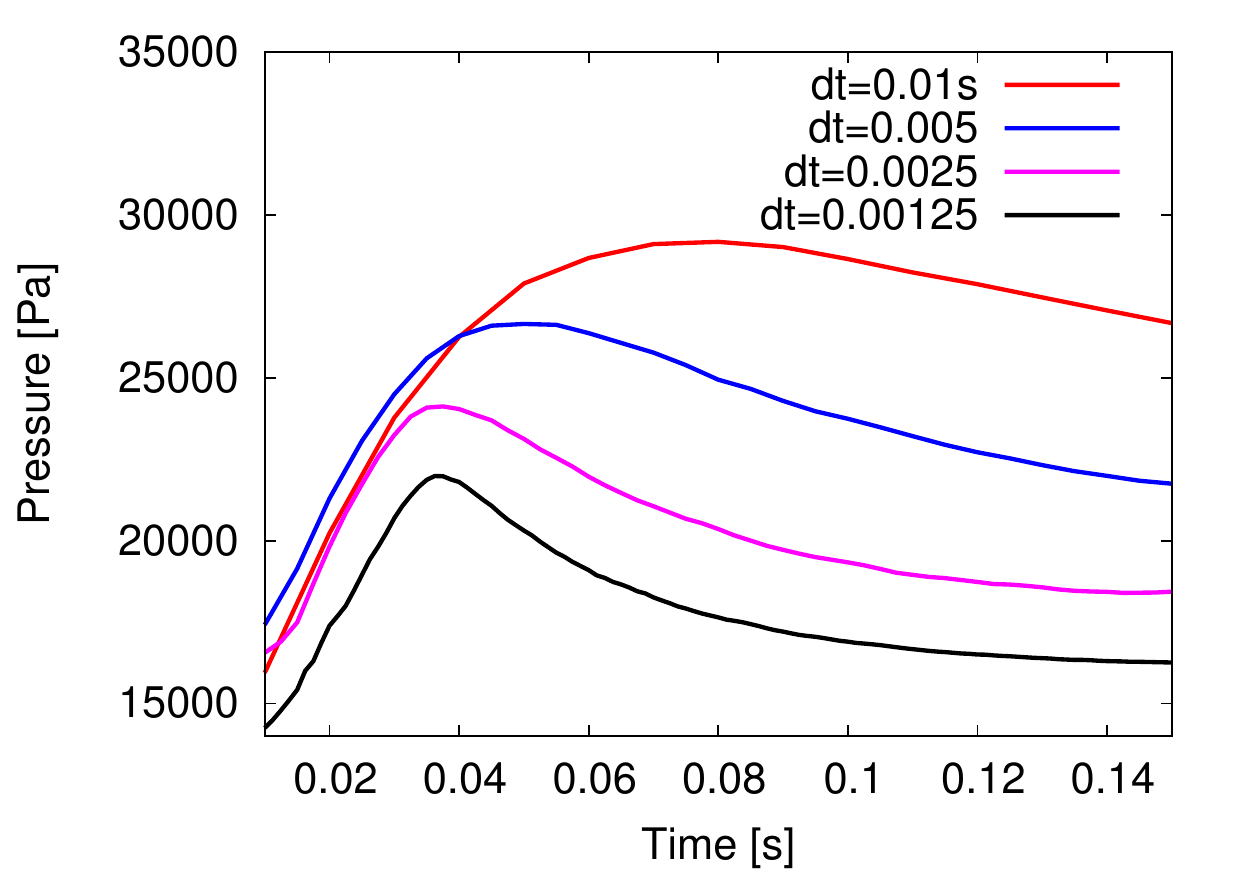}}
{\includegraphics[width=7cm]{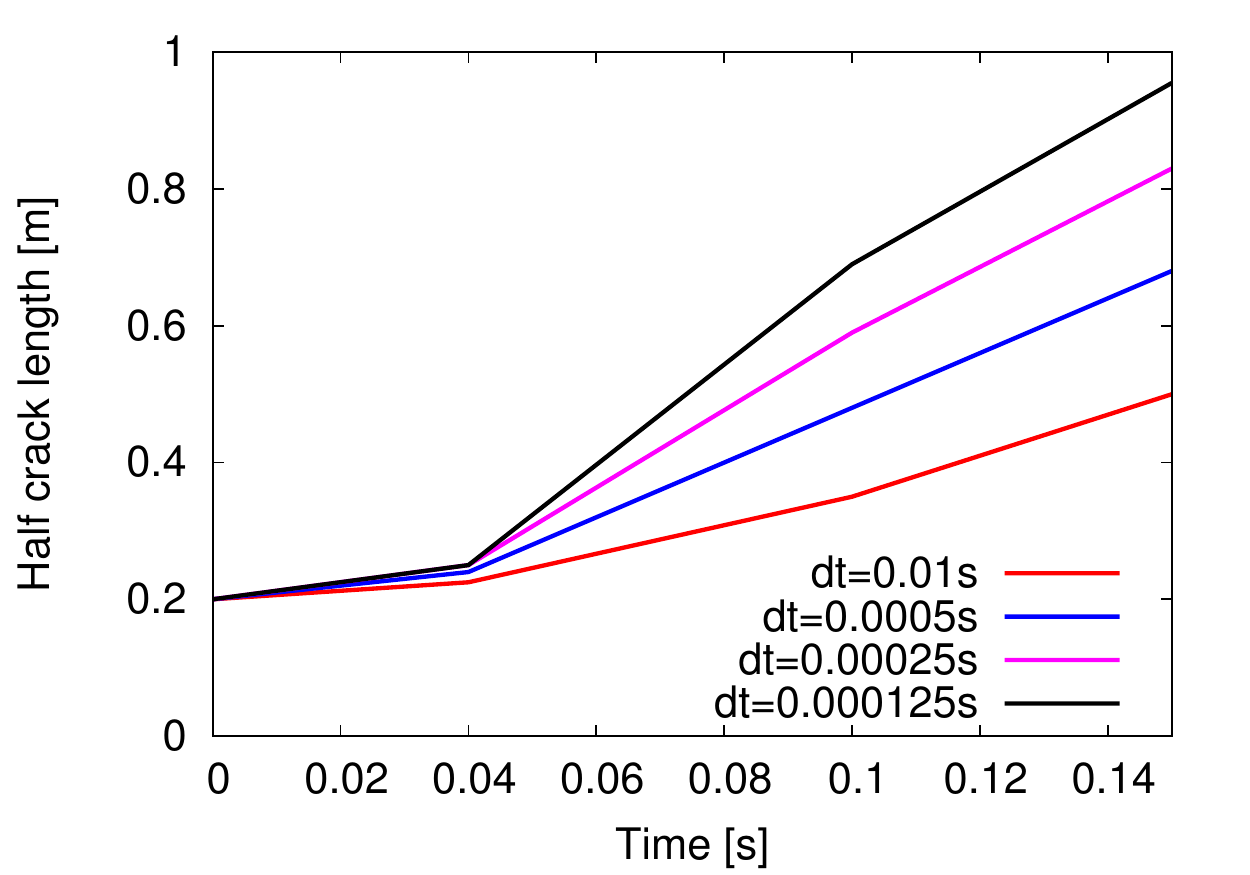}}
\caption{Example 3: Maximal pressure evolution and half crack length
 evolution on the finest mesh level. The  
time step sizes $\delta t$ are refined in order to study convergence 
in time. The maximal (theoretical) half length would be $2m$ (the boundary of the domain),
and the fractures stop growing towards $1.9m$.
Spatial refinement is not considered 
 since both $\eps$ and $h$ are varied via $\eps = 2h$ and convergence 
 to common values cannot be expected.
{ In the bottom row, zoom-ins are provided showing more clearly
  temporal convergence although it is very slow.}}
\label{ex_3_fig_2}
\end{figure}

\paragraph{Quantities of interest}
We study the following cases and goal functionals:
\begin{itemize}
\item the number of GMRES iterations within the Newton solver.
\item the number of fixed-stress iterations and 
 Newton solves for the displacement/phase-field system (Figure \ref{ex_3_fig_1});
\item the maximum pressure and 
the crack length/pattern evolution over time (Figure \ref{ex_3_fig_2}); 
\end{itemize}

\begin{figure}[H]
\centering
{\includegraphics[width=4cm]{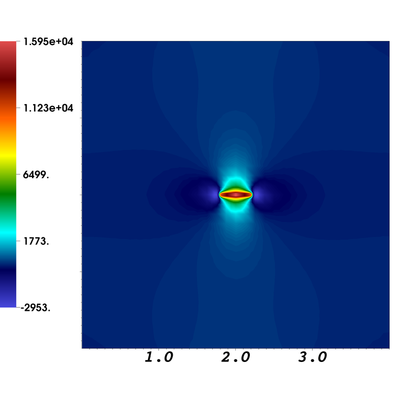}}
{\includegraphics[width=4cm]{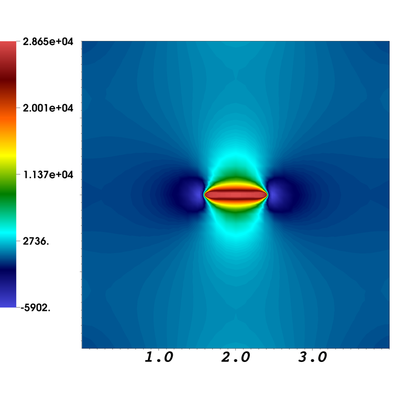}}
{\includegraphics[width=4cm]{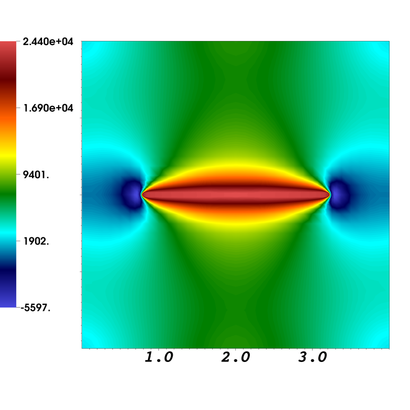}}
{\includegraphics[width=4cm]{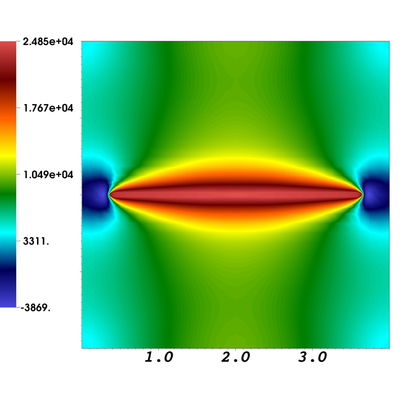}}\\
{\includegraphics[width=4cm]{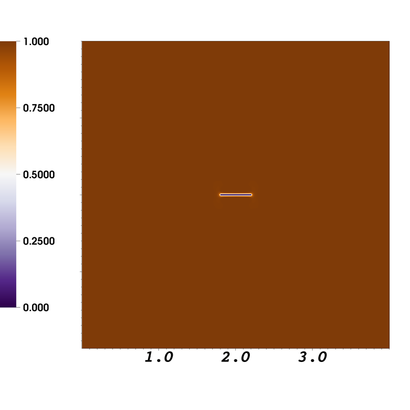}}
{\includegraphics[width=4cm]{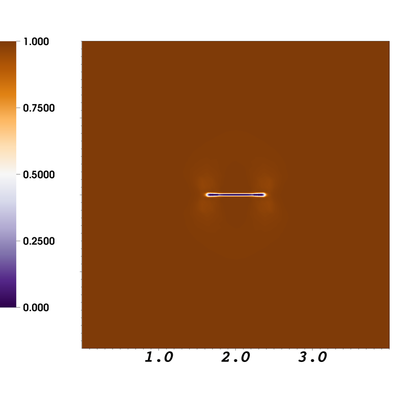}}
{\includegraphics[width=4cm]{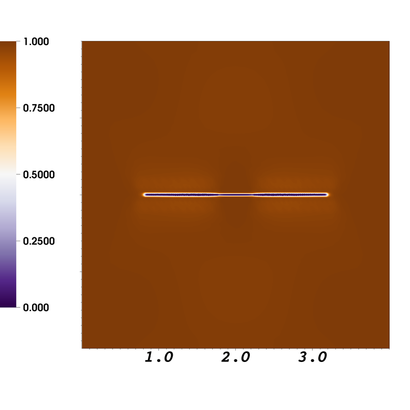}}
{\includegraphics[width=4cm]{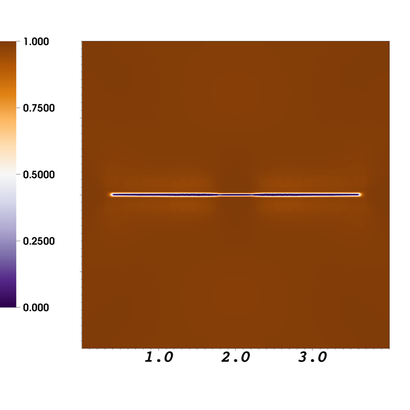}}\\
{\includegraphics[width=4cm]{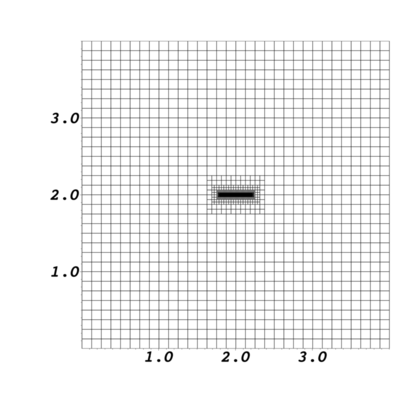}}
{\includegraphics[width=4cm]{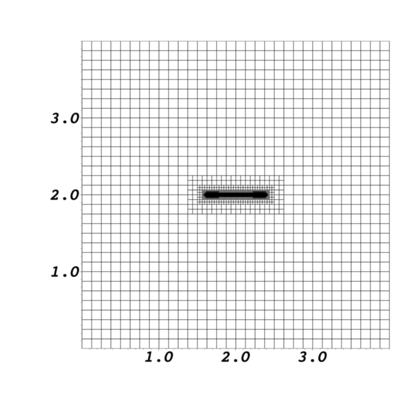}}
{\includegraphics[width=4cm]{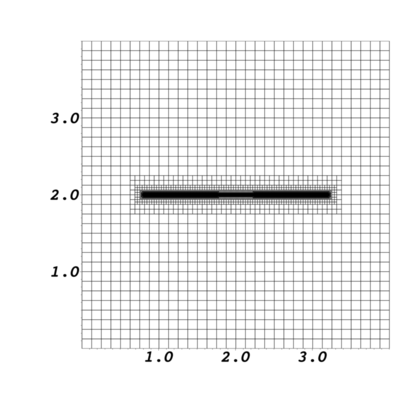}}
{\includegraphics[width=4cm]{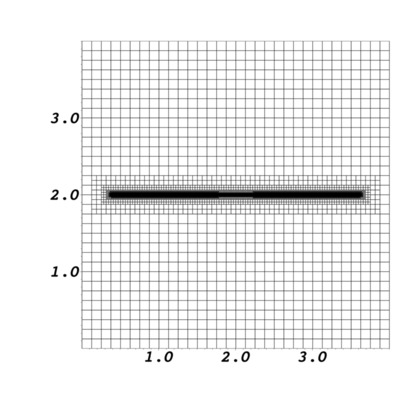}}
\caption{Example 3: Pressure evolution including negative pressure at the 
fracture tips (top), crack propagation (middle), 
and adaptive mesh evolution (bottom) at $T=0.01,0.1,0.4,0.6s$ and 
on the finest mesh level Ref. 5.}
\label{ex_3_fig_3}
\end{figure}

\paragraph{Discussion of findings}
The average number of Newton iterations for the displacement/phase-field 
system is $4-8$ iterations per mesh per time step. 
The average number of GMRES iterations is $10-50$ but can go up
in certain steps (just before the Newton tolerance is satisfied)
up to $100-150$, which is still acceptable.
To study fixed-stress iterations, computations are performed 
on different mesh levels in order to see the dependence of the number of 
mesh cells. Moreover, the presentation is divided into the number of 
fixed-stress iterations per mesh per time step and secondly, 
into the accumulated number (summing up all predictor-corrector mesh
refinements) per time step, see Figure \ref{ex_3_fig_1}. 
In Figure \ref{ex_3_fig_2}, we study temporal convergence 
of two quantities of interest; namely, the pressure and the 
fracture length. Here we identify convergence although it 
is slow. With regard to the complexity of the overall problem,
it is however a major accomplishment to obtain in the first place 
temporal convergence. This is an important finding with regard 
to the computational stability of the proposed framework.
{
We remark that almost identical computational results were observed by
using either Formulation \ref{form:level_set} or Formulation
\ref{form_level_set_by_PFF}. 
However, Formulation \ref{form:level_set} is more expensive
since an additional scalar-valued problem has to be solved.
}
Finally, in Figure \ref{ex_3_fig_3}, the pressure, the crack length 
in terms of the phase-field variable, and the locally refined mesh are
visualized.

\begin{figure}[H]
\centering
{\includegraphics[width=7cm]{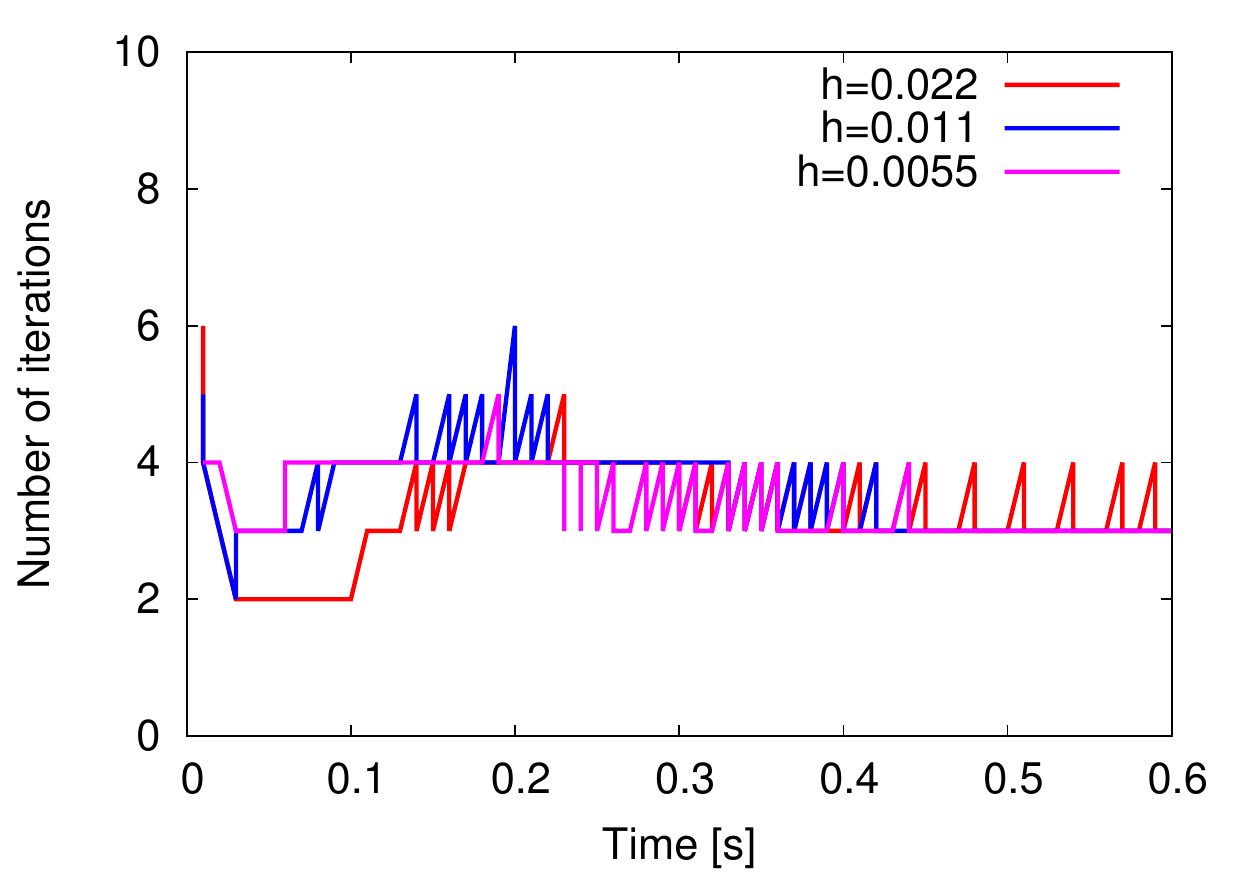}}
{\includegraphics[width=7cm]{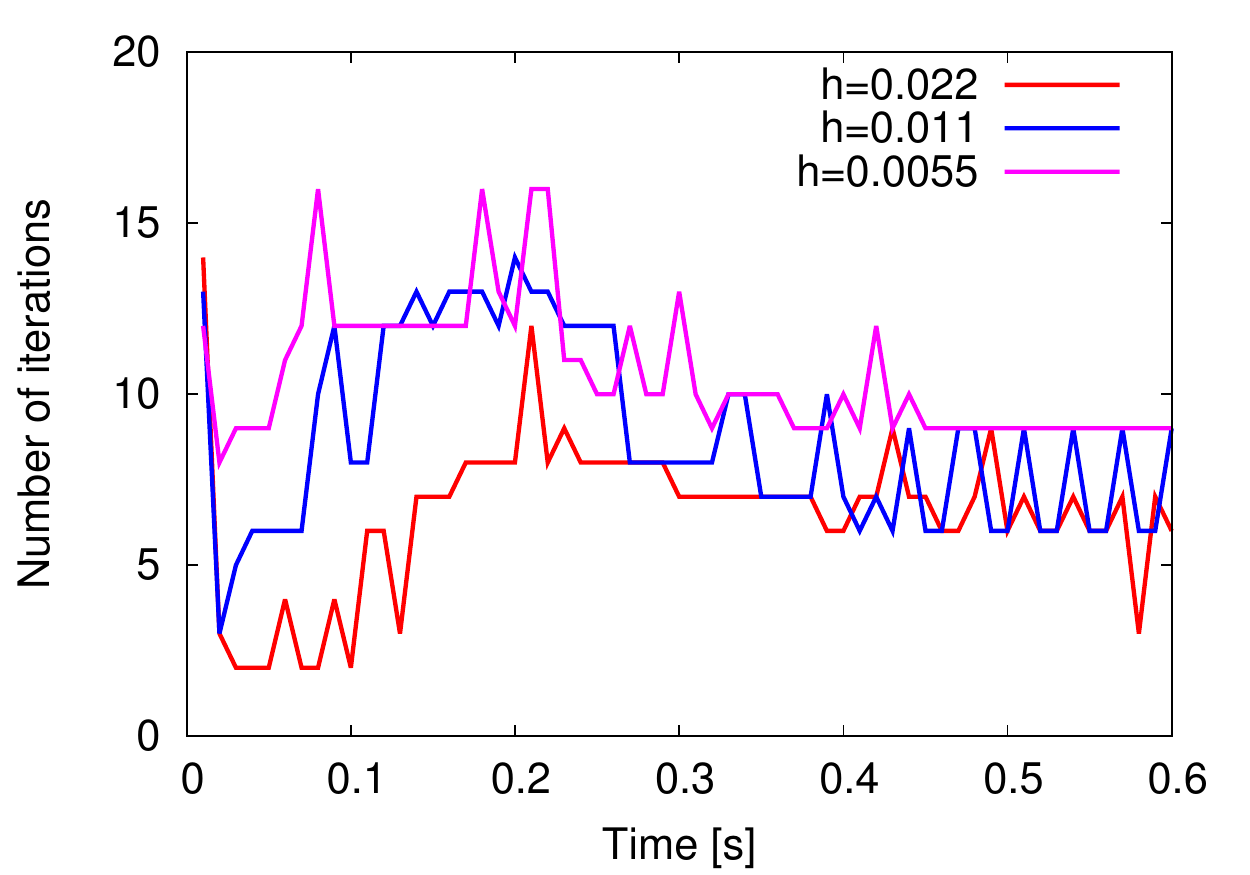}}
\caption{Example 3: Number of fixed-stress iterations per mesh per time step
  (at left). At right, the accumulated number of all predictor-corrector mesh
  levels per time step is displayed. The time step size in this test case $\delta t=0.01s$.}
\label{ex_3_fig_1}
\end{figure}

\subsection{Example 4: Fracture networks in homogeneous and heterogeneous porous media}
In this example, we study the fixed-stress algorithm for 
multiple fractures in homogeneous and heterogeneous porous media.
In total, we have three test cases: homogeneous, a heterogeneous 
example in which the Lam\'e parameters are varied, and a third 
example in which additionally the reservoir permeability is non-homogeneous.

\paragraph{Configuration}
We deal with the following geometric data:
$\Omega = (\SI{0}{\metre},\SI{10}{\metre})^2$ and three initial cracks.

\begin{figure}[H]
\centering
{\includegraphics[width=13cm]{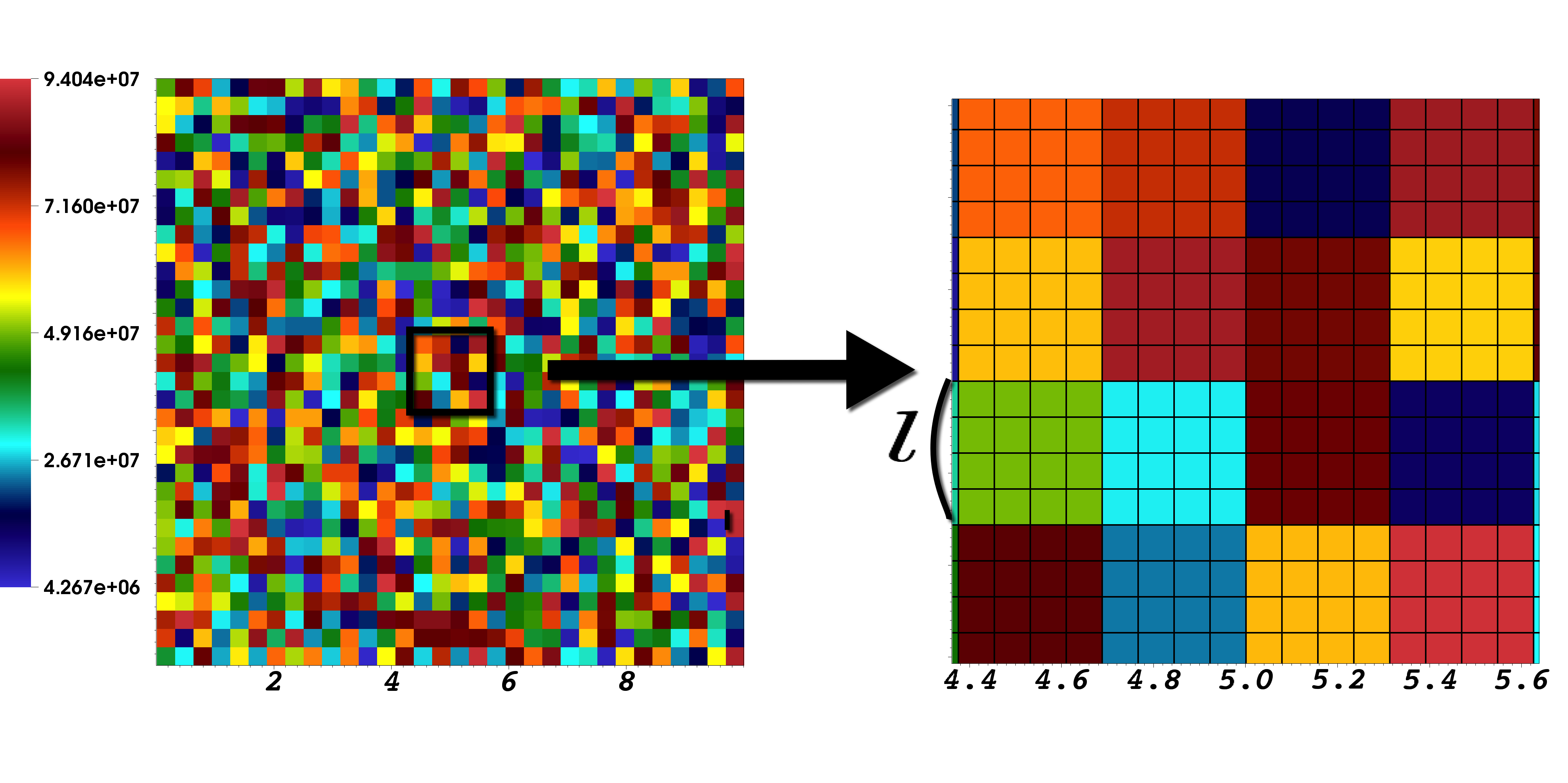}}
\caption{Example 4: Heterogeneous solid. Display of $\mu$. These data are synthetic using a random distribution. In particular, the heterogeneities are chosen as such that $\eps$ can resolve them, i.e., $l > \eps = 2h$, where $l$ is the length-scale of the material inhomogeneity (at right). 
}
\label{ex_4_fig_1}
\end{figure}

\paragraph{Parameters}
The fracture toughness is chosen as
$G_c = \SI{1}{\newton\per\metre}$. The mechanical parameters are
Young's modulus and Poisson's ratio $E = \SI{e+8}{\pascal}$ and $\nu_s = 0.2$
for the homogeneous case. In the heterogeneous setting,
we have $\SI{e+7}{\pascal}\leq E \leq \SI{e+8}{\pascal}$. These 
heterogeneities are chosen as such that the length-scale parameter
$\eps$ can resolve them; see Figure \ref{ex_4_fig_1}.

The regularization parameters are chosen as $\eps = 2h$
and $\kappa = \num{e-10}h$.
Biot's coefficient is $\alpha = 1$.
Furthermore $q_F = \SI{5}{m^3/s}$
and $M = \SI{1e+8}{Pa}, c_F = \SI{1e-8}{Pa}$. The viscosities are chosen as
$\nu_R = \nu_F = \SI{1e-3}{Ns/m^2}$. The reservoir permeability is $K_R = \SI{1}{d}$
in the homogeneous case and varies $0.1d\leq K_R\leq 1d$ in 
the heterogeneous setting.
and the density is $\rho_F^0 = \SI{1}{kg/m^3}$.
The time step size is $\delta t = \SI{1e-2}{s}$ and the final time is not 
specified and rather taken when all fractures joined. This 
event takes place between $0.25s\leq T\leq 0.3s$.
Furthermore, $TOL_{FS} = 10^{-4}$ (for the pressure and the displacements), 
whereas the phase-field tolerance is chosen as $TOL_{FS} = 10^{-2}$.
In fact the convergence of the phase-field variable is much harder 
for multiple fractures and heterogeneous media than in the previous 
examples.

\paragraph{Quantities of interest}
In this example, we observe the crack pattern, the pressure distribution, 
and fixed-stress iterations.

\begin{figure}[H]
\centering
{\includegraphics[width=4cm]{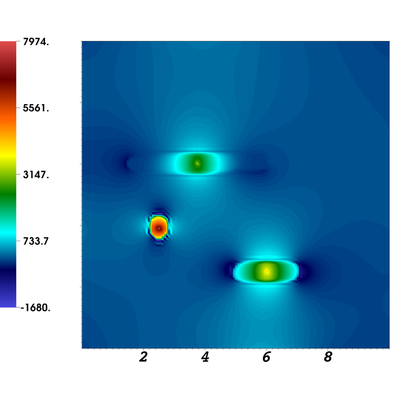}}
{\includegraphics[width=4cm]{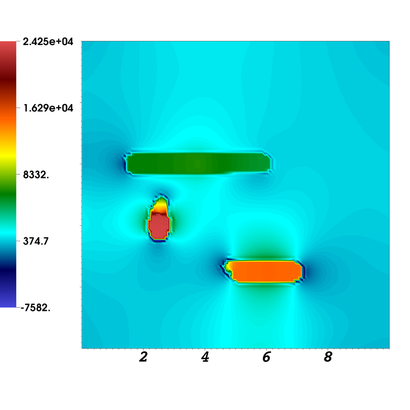}}
{\includegraphics[width=4cm]{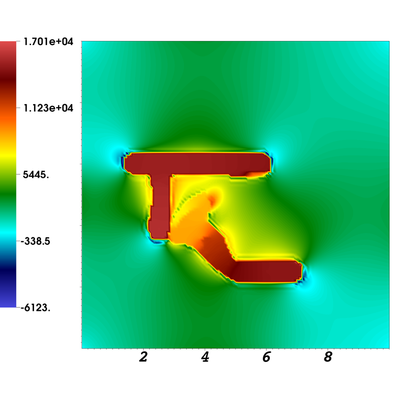}}
{\includegraphics[width=4cm]{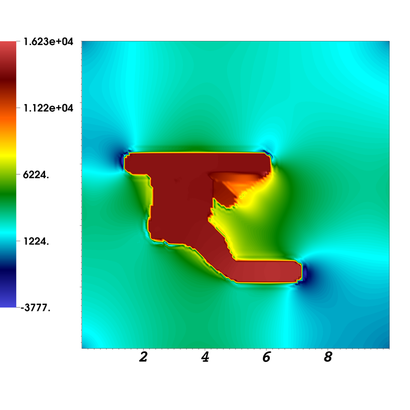}}\\
{\includegraphics[width=4cm]{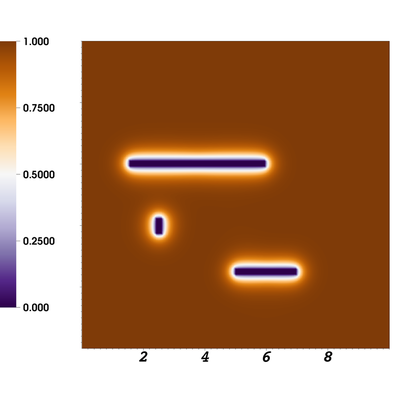}}
{\includegraphics[width=4cm]{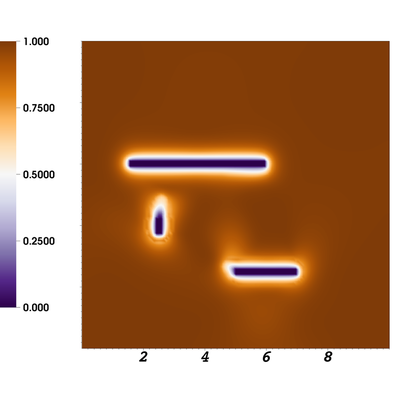}}
{\includegraphics[width=4cm]{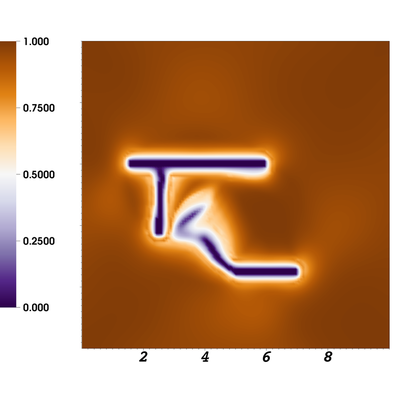}}
{\includegraphics[width=4cm]{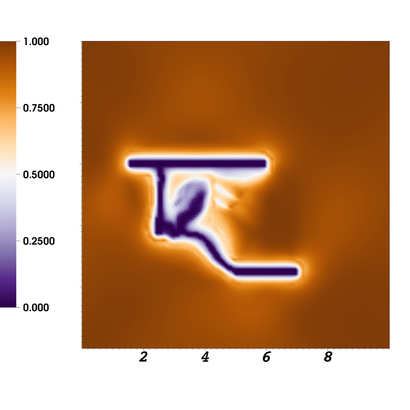}}
\caption{Example 4: Homogeneous test case. Display of pressure (top) and 
fracture pattern at $T=0s,0.1s,0.2s,0.25s$.}
\label{ex_4_fig_2}
\end{figure}

\begin{figure}[H]
\centering
{\includegraphics[width=4cm]{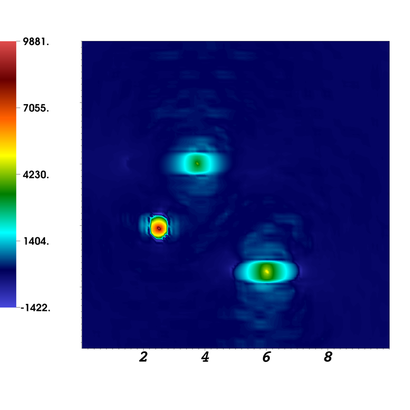}}
{\includegraphics[width=4cm]{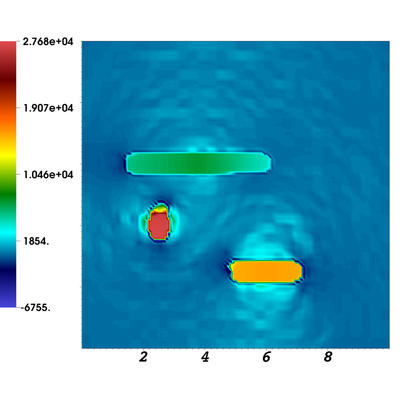}}
{\includegraphics[width=4cm]{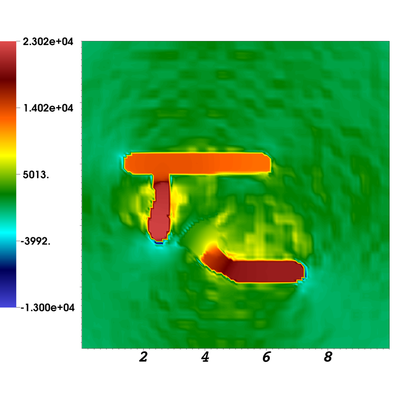}}
{\includegraphics[width=4cm]{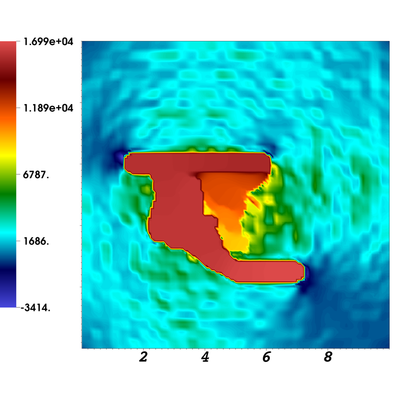}}\\
{\includegraphics[width=4cm]{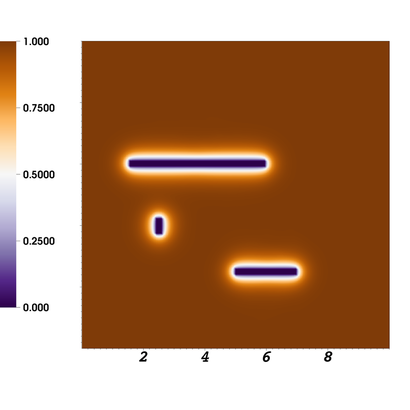}}
{\includegraphics[width=4cm]{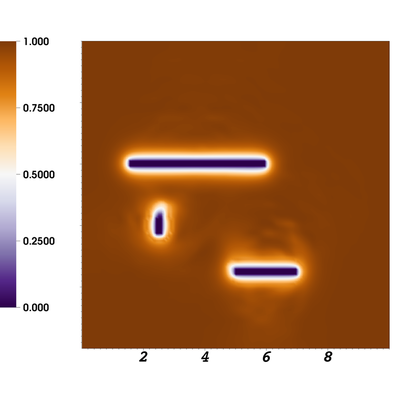}}
{\includegraphics[width=4cm]{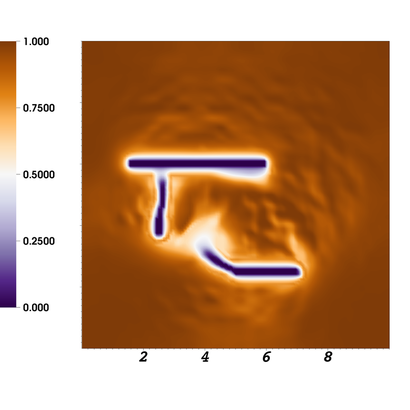}}
{\includegraphics[width=4cm]{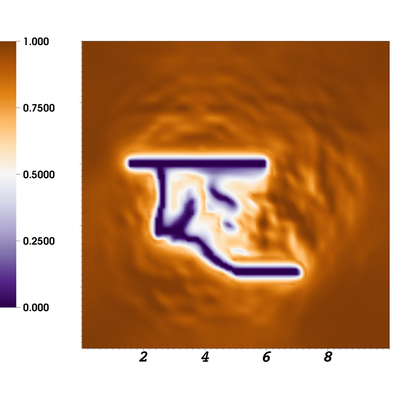}}
\caption{Example 4: Heterogeneous test case with varying Lam\'e parameters
and homogeneous reservoir permeability. Display of pressure (top) and 
fracture pattern at $T=0s,0.1s,0.2s,0.3s$.}
\label{ex_4_fig_3}
\end{figure}

\paragraph{Discussion of findings}
In Figure \ref{ex_3_fig_1}, the number of fixed-stress iterations 
per time step is shown.
The evolution of the pressure and the fracture patterns at different times
are displayed in Figure \ref{ex_4_fig_2} and Figure \ref{ex_4_fig_3}.
The average number of non-linear iterations of the 
semi-smooth Newton solver for the three test cases 
are $8-10$, with in average $10-20$ linear 
GMRES iterations. Here, we do not observe a significant 
difference between homogeneous and heterogeneous media.
{ Furthermore, we observe again the 
same crack pattern for both level-set formulations, but 
using Formulation \ref{form_level_set_by_PFF}  the final 
shape is reached earlier. Due to the complexity 
of this test (multiple fractures and heterogeneous materials)
further future investigations are definitely necessary.}

\begin{figure}[H]
\centering
{\includegraphics[width=8cm]{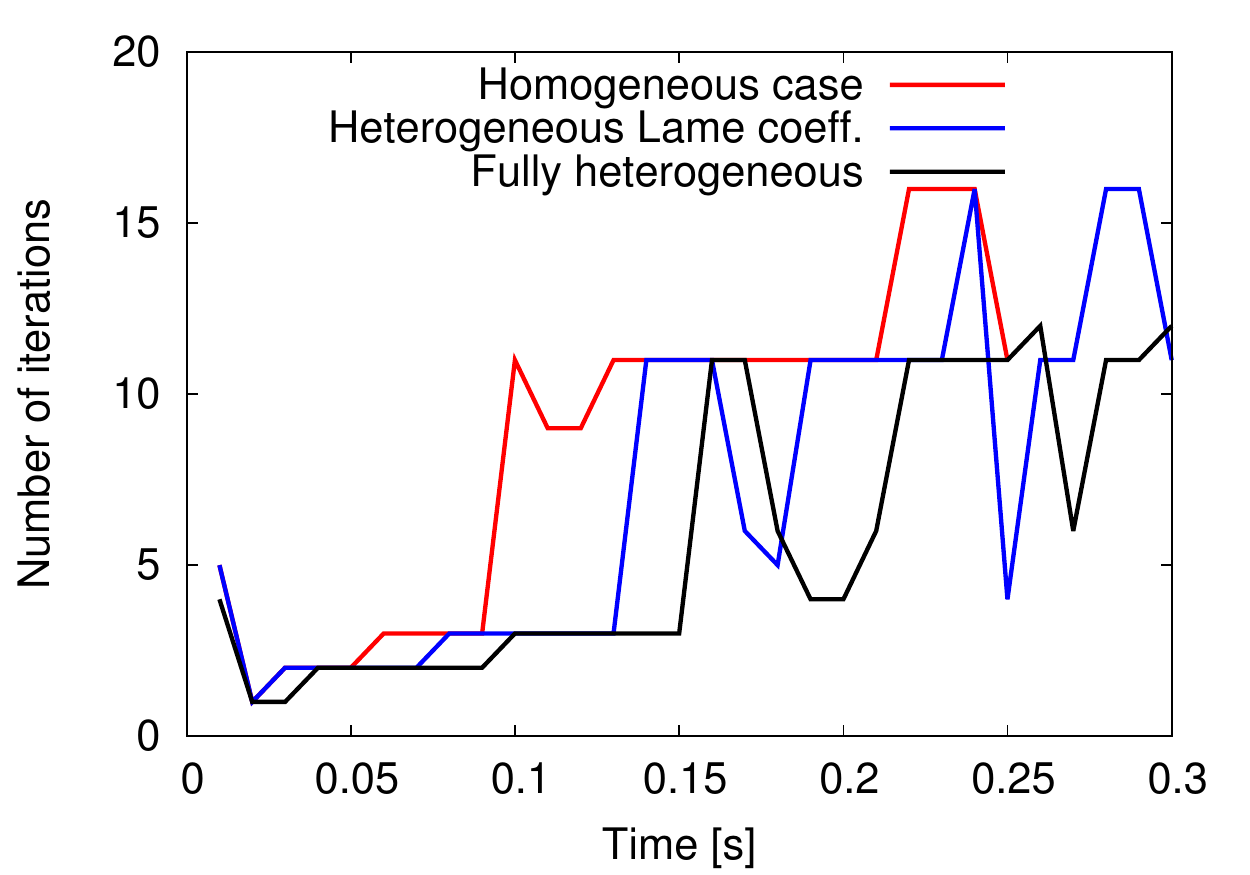}}
\caption{Example 4: The number of fixed-stress iterations per time step. 
The time step size in this test case $\delta t=0.01s$.}
\label{ex_3_fig_1}
\end{figure}

\begin{figure}[H]
\centering
{\includegraphics[width=4cm]{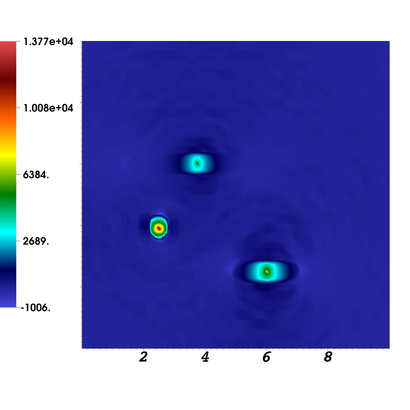}}
{\includegraphics[width=4cm]{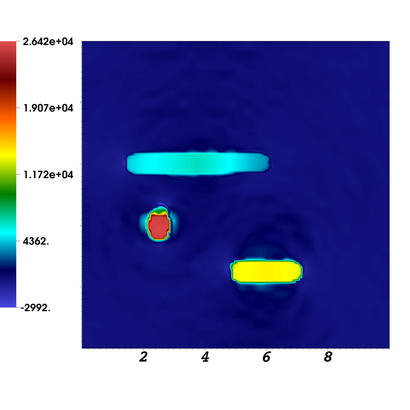}}
{\includegraphics[width=4cm]{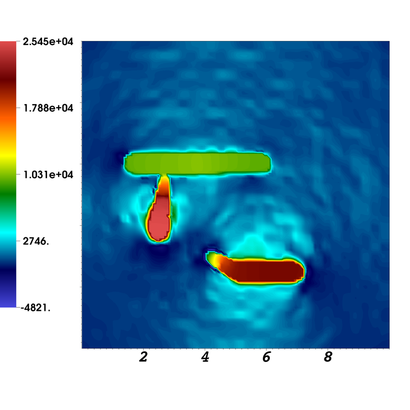}}
{\includegraphics[width=4cm]{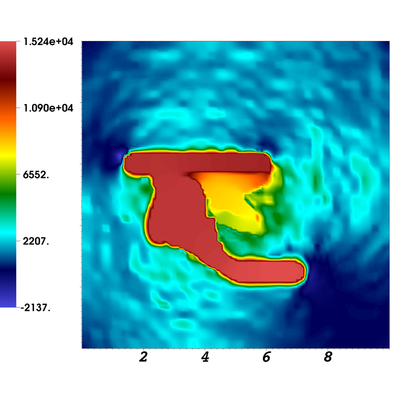}}\\
{\includegraphics[width=4cm]{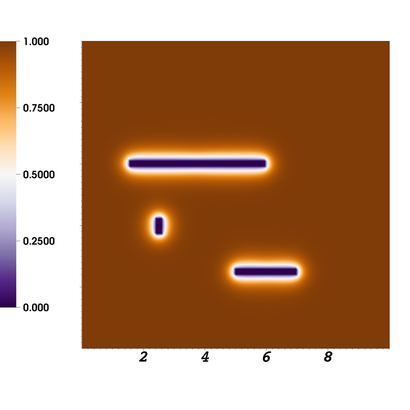}}
{\includegraphics[width=4cm]{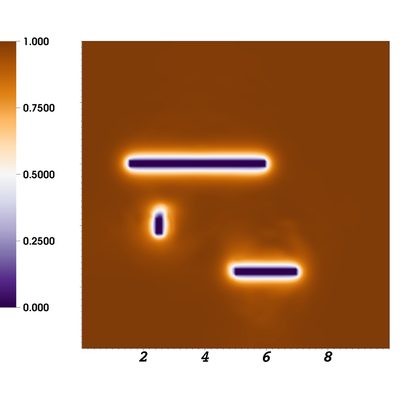}}
{\includegraphics[width=4cm]{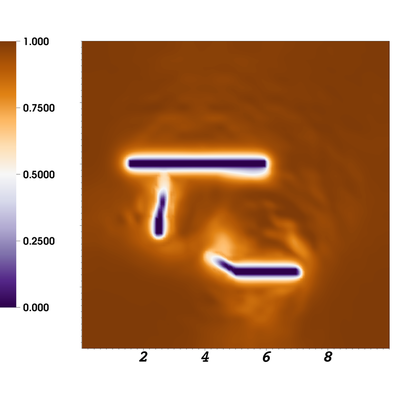}}
{\includegraphics[width=4cm]{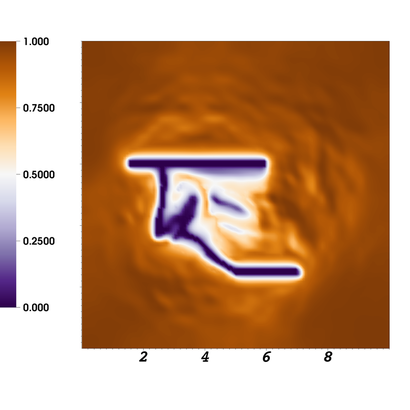}}
\caption{Example 4: Fully heterogeneous test case. Display of pressure (top) and 
fracture pattern at $T=0s,0.1s,0.2s,0.3s$.}
\label{ex_4_fig_4}
\end{figure}

\subsection{Example 5: Propagating penny-shaped fracture in 3D}
\label{sec:ex_5}
\begin{figure}[!h]
\centering
\begin{subfigure}[b]{0.225\textwidth}
\includegraphics[width=\textwidth]{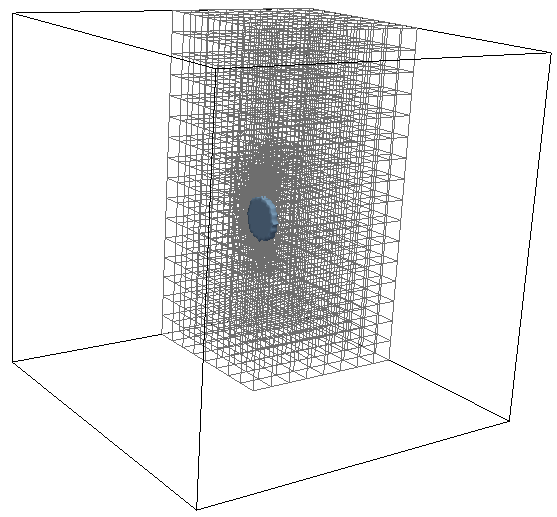}
\caption{Setup}
\label{fig:num:sneddon_set_a}
\end{subfigure}
\begin{subfigure}[b]{0.225\textwidth}
\includegraphics[width=\textwidth]{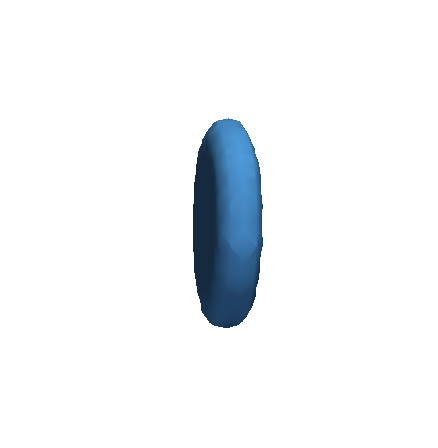}
\caption{$n=30$}
\label{fig:num:sneddon_set_b}
\end{subfigure}
\begin{subfigure}[b]{0.225\textwidth}
\includegraphics[width=\textwidth]{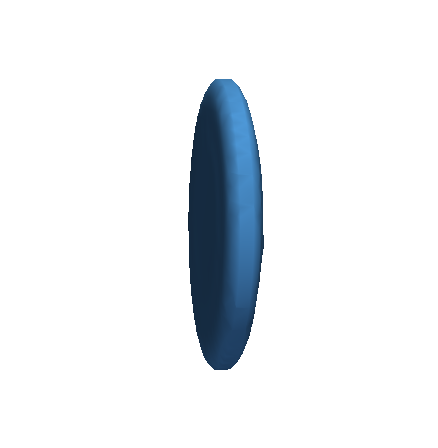}
\caption{$n=50$}
\label{fig:num:sneddon_set_c}
\end{subfigure}
\begin{subfigure}[b]{0.225\textwidth}
\includegraphics[width=\textwidth]{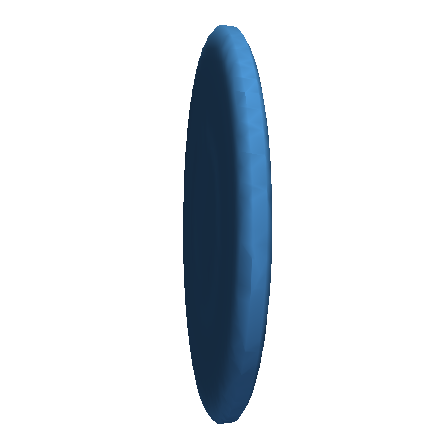}
\caption{$n=100$}
\label{fig:num:sneddon_set_d}
\end{subfigure}

\caption{Example \ref{sec:ex_5}: (a) The initial penny shape crack is centered at $(\SI{5}{\metre},\SI{5}{\metre},\SI{5}{\metre})$ on $y=\SI{5}{\metre}-$plane with the mesh refinement. (b)-(d) illustrate the propagating fracture for each time step.}
\label{fig:num:sneddon_set}
\end{figure}

In this final example, we consider a penny shaped fracture in a three dimensional domain $\Omega = (0,\SI{4}{\metre})^3$. 
The horizontal initial penny shape crack is centered at $(\SI{2}{\metre},\SI{2}{\metre},\SI{2}{\metre})$ on $y=\SI{2}{\metre}-$ plane and we refine around the crack; see Figure \ref{fig:num:sneddon_set_a} for the setup.
Initial and boundary conditions are same as previous examples and here the physical parameters are given.
The fracture toughness is chosen as
$G_c = \SI{1}{\newton\per\metre}$,
Young's modulus and Poisson's ratio as $E = \SI{e+8}{\pascal}$ and $\nu_s = 0.2$, respectively. 
The regularization parameters are chosen as $\eps = 2h_{\min}$
and $\kappa = \num{e-10}h$.
The Biot coefficient and Biot's modulus are set to $\alpha = 1$ 
and $M = \SI{1e+8}{Pa}$, respectively.
Furthermore we assume a slightly incompressible fluid with 
$c_F = \SI{1e-8}{Pa}$ and the viscosities are chosen as
$\eta_R = \eta_F = \SI{1e-3}{Ns/m^2}$ with the injection rate $q_F =\SI{2}{m^3/s}$. 
The reservoir permeability is $K_R = \SI{1}{d}$ 
and the density is $\rho_F^0 = \SI{1}{kg/m^3}$.
Furthermore,
$h_{\min}= \SI{0.05}{\metre}$, $\delta t = \SI{1e-2}{s}$, and 
 $TOL_{FS} = 10^{-3}$.

\begin{figure}[!h]
\centering
\begin{subfigure}[b]{0.4\textwidth}
\includegraphics[width=\textwidth]{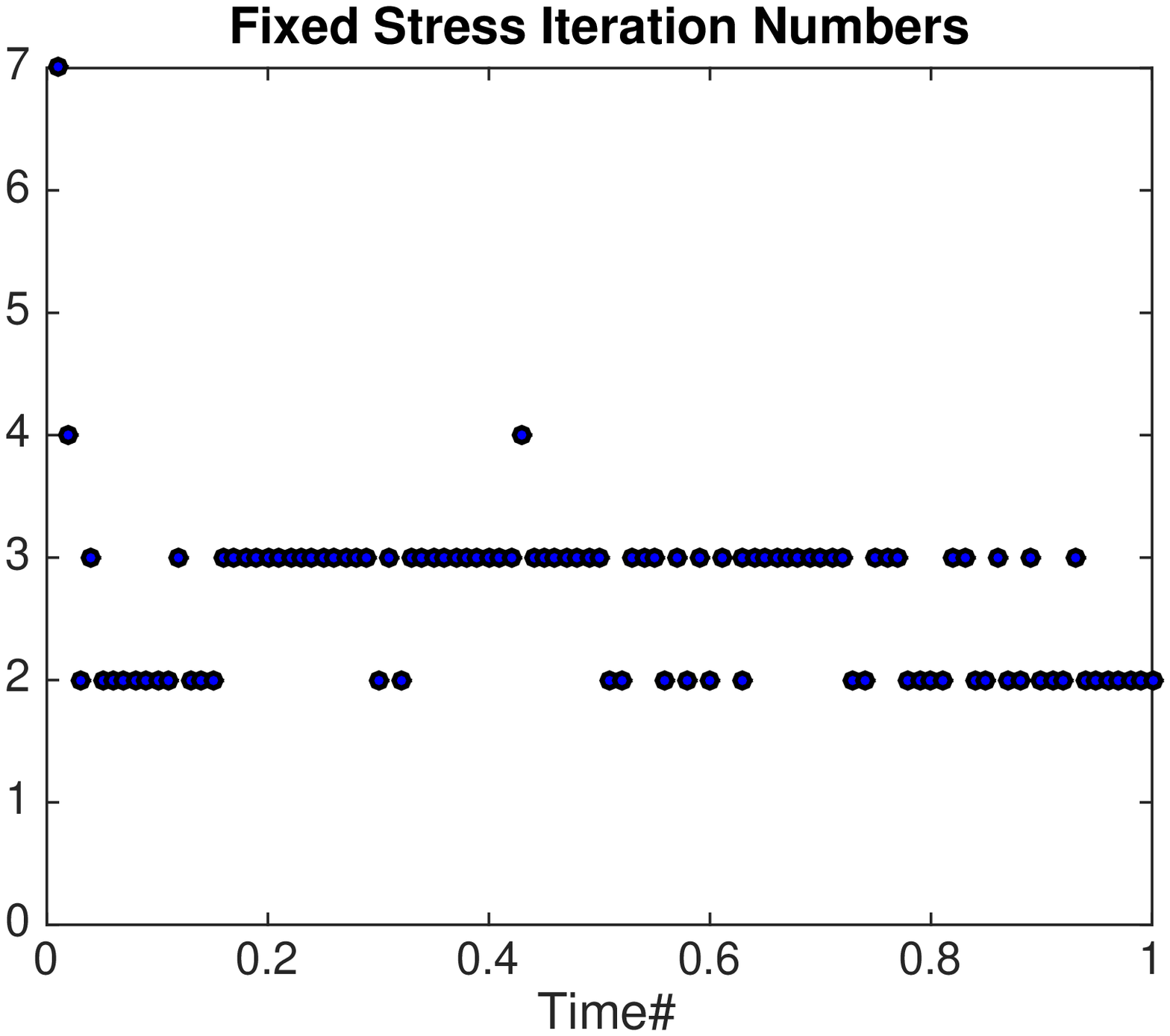}
\caption{Fixed stress iteration number}
\end{subfigure}
\begin{subfigure}[b]{0.4\textwidth}
\includegraphics[width=\textwidth]{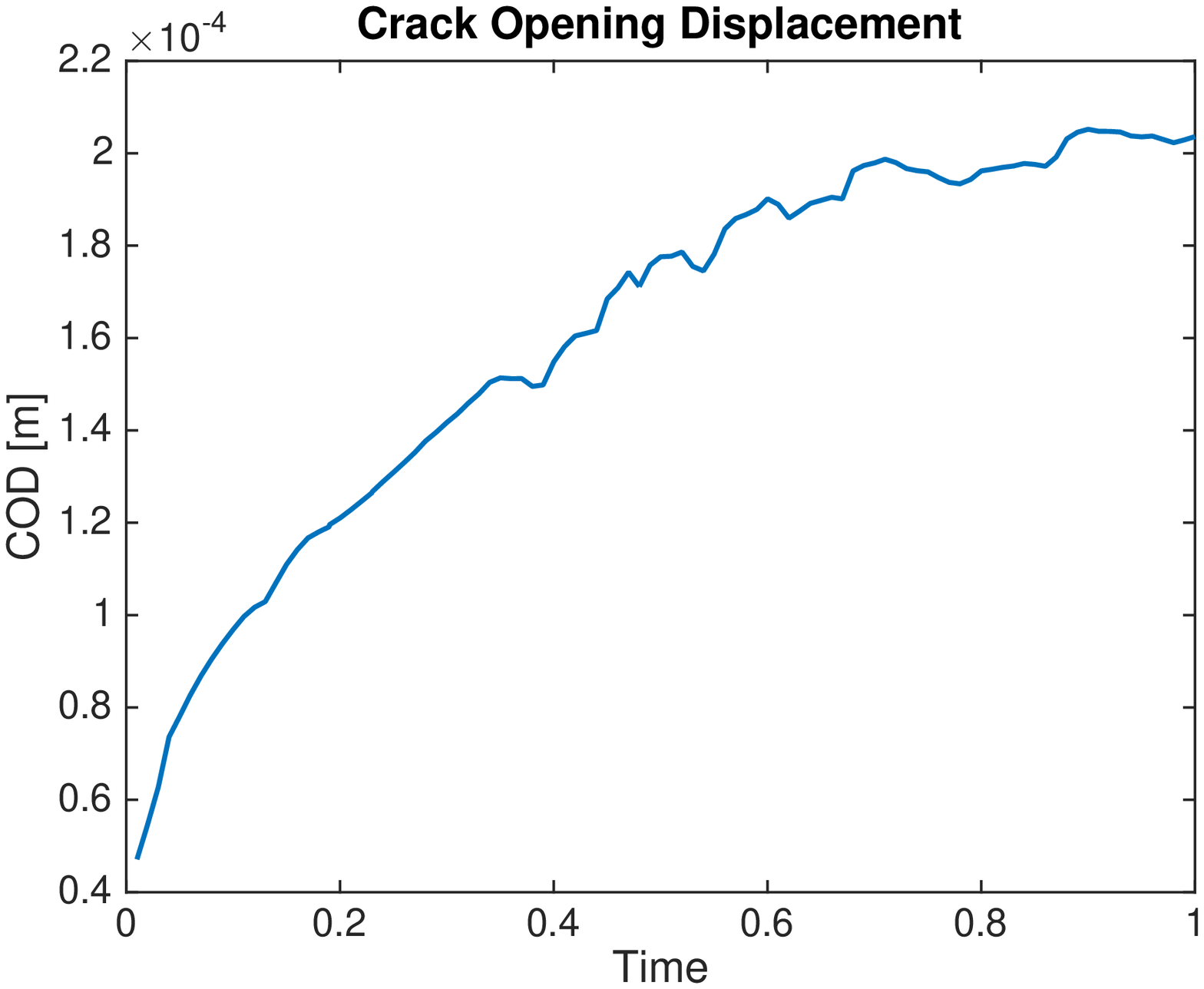}
\caption{Crack opening displacement $\|w \|_{L^{\infty}(\Lambda)}$.}
\end{subfigure}
\caption{(a) Fixed stress iteration number over the time with a three
  dimensional MPI computation is shown. (b) crack width opening over the time
  for the propagating fracture.}
\label{fig:ex3D_FS}
\end{figure}

Figure \ref{fig:num:sneddon_set_b}-\ref{fig:num:sneddon_set_d} illustrate the propagating fracture for each time. 
The crack opening displacement and the fixed stress iteration number over the time for  propagating fracture are shown in Figure \ref{fig:ex3D_FS}.
{
We note that the almost identical computational results were observed by
employing either Formulation \ref{form:level_set} or Formulation
\ref{form_level_set_by_PFF} to compute the level-set as shown in previous
example. However, Formulation \ref{form_level_set_by_PFF} is
computationally more  efficient especially in three dimensional cases.   }

\section{Conclusions}
In this paper, we presented fixed-stress splitting for 
fractured porous media using a phase-field technique. 
Several examples were consulted in order to show the performance 
of the algorithmic techniques. Despite the complexity of 
the entire problem, we could observe spatial and temporal 
convergence of selected quantities of interest. This is 
a major step towards the computational stability and reliability 
of our proposed method.
Moreover, we could observe typical properties of 
fluid-filled fractures, namely a negative pressure at the 
fracture tips, which are not present when Biot's coefficient 
is zero.
Moreover, we investigated the solver iteration numbers for 
the linear iterative GMRES solver, the nonlinear Newton solver 
of the displacement/phase-field system, and the fixed stress iterations
between flow and mechanics. 
For our settings, we obtained efficient iterations numbers as shown
in our numerical examples.
However, in Example 4, heterogeneous materials, we 
observed that the convergence is dominated by the phase-field variable whereas the pressure and the displacements converges well.
We finally mention that the level-set width computation shows a novel way to obtain accurate width values inside the fracture region.
 This holds in particular true for homogeneous test cases (Examples 1-3 and
 Example 5). For heterogeneous 
tests and multiple fractures (Example 4), we also obtained good results but it is still an open question whether  the methodology works for arbitrary heterogeneous materials which goes  beyond the current paper and is left for future research.

\section*{Acknowledgments}
The authors want to thank Brice Lecampion, Emmanuel Detournay,  
Alf Birger Rustad, 
H\r{a}kon H\o gst\o l,
and Ali Dogru
for providing information and discussions on 
fluid-filled fractures and the resulting pressure behavior.
The research by S. Lee and  M. F. Wheeler was partially supported by 
{a} DOE grant DE-FG02-04ER25617,
{a} Statoil grant STNO-4502931834, and 
{an} Aramco grant UTA 11-000320. 
T. Wick would like to thank  
the JT Oden Program of the {Institute for Computational  Engineering and Science  (ICES)} and the Center for Subsurface Modeling (CSM), UT Austin for funding and hospitality during his visit in April  2016.


\bibliographystyle{abbrv}
\bibliography{lit}


\end{document}